\documentclass[11pt]{article}

\usepackage{amsmath, amssymb, amsthm, graphicx}
\usepackage[active]{srcltx}

\title{Rainbow Hamilton cycles in random graphs}

\author{
Alan Frieze \thanks{Department of Mathematical
Sciences, Carnegie Mellon University, Pittsburgh, PA 15213, email:
alan@random.math.cmu.edu. Research supported in part by NSF award
DMS-0753472.}
\and
Po-Shen Loh
\thanks{Department of Mathematical Sciences, Carnegie Mellon University,
Pittsburgh, PA 15213, e-mail: ploh@cmu.edu.}
}

\date{}

\newtheorem{theorem}{Theorem}[section]
\newtheorem*{theorem*}{Theorem}
\newtheorem{definition}[theorem]{Definition}

\newtheorem{lemma}[theorem]{Lemma}
\newtheorem{corollary}[theorem]{Corollary}
\newtheorem{question}[theorem]{Question}

\newcommand{\proofstart}{\noindent \textbf{Proof.}\, }
\newcommand{\proofend}{\hfill $\Box$}

\newcommand{\set}[1]{\left\{ #1 \right\}}
\newcommand{\bfrac}[2]{\left(\frac{#1}{#2}\right)}
\newcommand{\ex}[1]{\exp\left\{ #1 \right\}}

\newcommand{\pr}[1]{\mathbb{P}\left[ #1 \right]}
\newcommand{\E}[1]{\mathbb{E}\left[ #1 \right]}

\newcommand{\bin}[1]{\text{\rm Bin}\left( #1 \right)}

\def\e{\epsilon}    \def\g{\gamma}
  \def\k{\kappa}
 \def\th{\theta}

\newcommand{\cE}{\mathcal{E}}

\newcommand{\whp}{\textbf{whp}}
\newcommand{\Gnp}{G_{n,p}}
\newcommand{\Dnp}{D_{n,p}}

\newcommand{\Gnpk}{G_{n,p,\kappa}}

\newcommand{\dist}{\text{dist}}
\newcommand{\rt}[1]{\overrightarrow{#1}}

\begin{document}
\maketitle

\begin{abstract}
  One of the most famous results in the theory of random graphs
  establishes that the threshold for Hamiltonicity in the
  Erd\H{o}s-R\'enyi random graph $G_{n,p}$ is around $p \sim
  \frac{\log n + \log \log n}{n}$.  Much research has been done to
  extend this to increasingly challenging random structures.  In
  particular, a recent result by Frieze determined the asymptotic
  threshold for a loose Hamilton cycle in the random 3-uniform
  hypergraph by connecting 3-uniform hypergraphs to edge-colored
  graphs.

  In this work, we consider that setting of edge-colored graphs, and
  prove a result which achieves the best possible first order
  constant.  Specifically, when the edges of $G_{n,p}$ are randomly
  colored from a set of $(1 + o(1)) n$ colors, with $p = \frac{(1 +
    o(1)) \log n}{n}$, we show that one can almost always find a
  Hamilton cycle which has the further property that all edges are
  distinctly colored (rainbow).
\end{abstract}

\section{Introduction}

Hamilton cycles occupy a position of central importance in graph
theory, and are the subject of countless results.  In the context of
random structures, much research has been done on many aspects of
Hamiltonicity, in a variety of random structures.  See, e.g., any of
\cite{Bolbook, Bol-hitting-time, BolFrieze-hitting-time, KS-Hamthresh,
RW-3reg} concerning Erd\H{o}s-R\'enyi random graphs and
random regular graphs, any of \cite{CF-2-inout,  DFrieze, McD-digraphs1,
  McD-digraphs2} regarding directed graphs, or any of the recent
developments \cite{DudekFr, FrH, FK-pack, FKL-pack} on uniform
hypergraphs.  In this paper we consider the existence of rainbow
Hamilton cycles in edge-colored graphs.  (A set $S$ of edges is called
\emph{rainbow}\/ if each edge of $S$ has a different color.)  There
are two general types of results in this area: existence \whp
\footnote{A sequence of events $\cE_n$ is said to occur {\em with high
    probability} (\whp) if $\lim_{n\to\infty} \Pr{\cE_n}=1$.}  under
random coloring and guaranteed existence under adversarial coloring.

When considering adversarial (worst-case) coloring, the guaranteed
existence of a rainbow structure is called an \emph{Anti-Ramsey}\/
property.  Erd\H{o}s, Ne\v{s}et\v{r}il, and R\"odl \cite{ENR}, Hahn
and Thomassen \cite{HT} and Albert, Frieze, and Reed \cite{AFR}
(correction in Rue \cite{Rue}) considered colorings of the edges of
the complete graph $K_n$ where no color is used more than $k$
times. It was shown in \cite{AFR} that if $k \leq n/64$, then there
must be a multi-colored Hamilton cycle.  Cooper and Frieze \cite{CF1}
proved a random graph threshold for this property to hold in almost
every graph in the space studied.

There is also a history of work on random coloring (see, e.g., any of
\cite{CF1, CF, FM, JW}), and it has recently become apparent that this
random setting may be of substantial utility.  Indeed, a result of
Janson and Wormald \cite{JW} on rainbow Hamilton cycles in randomly
edge-colored random regular graphs played a central role in the recent
determination of the threshold for loose Hamiltonicity in random
3-uniform hypergraphs by Frieze \cite{FrH}.  Roughly speaking, a
hyperedge (triple of vertices) can be encoded by an ordinary edge
(pair of vertices), together with a color.  Hence, a random 3-uniform
hypergraph gives rise to a randomly edge-colored random graph.  We
will discuss this further in Section \ref{sec:graphs-to-hypergraphs},
when we use the reverse connection to realize one part of our new
result.

Let us now focus on the random coloring situation, where we consider
the following model.  Let $\Gnpk$ denote a randomly colored random
graph, constructed on the vertex set $[n]$ by taking each edge
independently with probability $p$, and then independently coloring it
with a random color from the set $[\k]$.  We are interested in
conditions on $n,p,\k$ which imply that $\Gnpk$ contains a rainbow
Hamilton cycle \whp.  The starting point for our present work is the
following theorem of Cooper and Frieze.

\begin{theorem*}
  (See \cite{CF}, Theorem 1.1.)  There exist constants $K_1$ and $K_2$
  such that if $p > \frac{K_1 \log n}{n}$ and $\k > K_2 n$, then
  $\Gnpk$ contains a rainbow Hamilton cycle \whp.
\end{theorem*}

\noindent The aim of this paper is to substantially strengthen the
above result by proving:

\begin{theorem}\label{th1}\
\begin{description}

\item[(a)] There exists a constant $K$ such that if $p > \frac{K \log
    n}{n}$, then for even $n$, $G_{n,p,n}$ contains a rainbow Hamilton
  cycle \whp.

\item[(b)] If $p = \frac{(1 + \e) \log n}{n}$ and $\k = (1+\th)n$,
  where $\e, \th > \frac{100}{\sqrt{\log \log n}}$, then $\Gnpk$
  contains a rainbow Hamilton cycle \whp.

\end{description}
\end{theorem}

To discuss the tightness of our main theorem, let us recall the
threshold for Hamiltonicity in $G_{n,p}$, established by
Koml\'os and Szemer\'edi \cite{KS-Hamthresh}.  We find that we must have $p > \frac{\log n + \log
  \log n + \omega(n)}{n}$ with $\omega(n) \rightarrow \infty$, or else
the underlying uncolored $G_{n,p}$ will not even be Hamiltonian.  We
also need at least $n$ colors to appear on the edges in order to have
enough colors for a rainbow Hamilton cycle.  Note that the earlier
result came within a constant factor of both of these minimum
requirements, while part (a) above achieves the absolute best possible
constraint on the number of colors, while still staying within a
constant factor of the minimally required number of edges (albeit only
for even $n$).

Part (b) drives both constants down to be best possible up to first
order, and for all values of $n$, regardless of parity.  We permit our
error terms $\e$ and $\th$ to decrease slowly, although we do not
expect our constraints on them to be optimal.  Our discussion above
shows that the trivial lower bound for $\e$ is around $\frac{\log \log
  n}{\log n}$.  Then, if $p \sim \frac{\log n}{n}$, we need at least
$n + \Omega(n^{1/2})$ colors just to ensure that \whp\ at least $n$
distinct colors occur on the $m \sim \frac{1}{2} n \log n$ edges in
the graph; hence, the trivial lower bound for $\th$ is around
$\frac{1}{\sqrt{n}}$.  We leave further exploration to future work,
and highlight a potential answer in our conclusion.

This paper is organized as follows.  We begin by establishing part (a)
of Theorem \ref{th1} in the next section.  Section 3 outlines the
proof of part (b), which is the main contribution of this paper.  The
proofs of the main steps follow in the section thereafter.  We
conclude in Section 5 with some remarks and open problems.  The
following (standard) asymptotic notation will be utilized extensively.
For two functions $f(n)$ and $g(n)$, we write $f(n) = o(g(n))$, $g(n)
= \omega(f(n))$, or $f(n) \ll g(n)$ if $\lim_{n \rightarrow \infty}
f(n)/g(n) = 0$, and $f(n) = O(g(n))$ or $g(n) = \Omega(f(n))$ if there
exists a constant $M$ such that $|f(n)| \leq M|g(n)|$ for all
sufficiently large $n$.  We also write $f(n) = \Theta(g(n))$ if both
$f(n) = O(g(n))$ and $f(n) = \Omega(g(n))$ are satisfied.  All
logarithms will be in base $e \approx 2.718$.

\section{Colored graphs and 3-uniform hypergraphs}
\label{sec:graphs-to-hypergraphs}

In this section we prove part (a) of our main theorem, and demonstrate
the connection between rainbow Hamilton cycles in graphs and loose
Hamilton cycles in 3-uniform hypergraphs.  Indeed, this will allow us
to realize part (a) as essentially a reformulation of the following
recent result of Frieze \cite{FrH}.  Let $H_{n,p;3}$ denote the random
3-uniform hypergraph where each potential hyperedge appears
independently with probability $p$.  In this object, a \emph{loose
  Hamilton cycle}\/ is a permutation of the vertices $(v_1, \ldots,
v_n)$ such that $\{v_1, v_2, v_3\}$, $\{v_3, v_4, v_5\}$, \ldots,
$\{v_{n-1}, v_n, v_1\}$ all appear as hyperedges ($n$ must be even).

\begin{theorem}
  (See \cite{FrH}, Theorem 1.)
  \label{th:loose-3}
  There is an absolute constant $K$ such that if $p > \frac{K \log
    n}{n^2}$, then for all $n$ divisible by 4, the random hypergraph
  $H_{n,p;3}$ on the vertex set $\{x_1, \ldots, x_{n/2}, y_1, \ldots,
  y_{n/2}\}$ contains a loose Hamilton cycle \whp.  Furthermore, one
  can find such a cycle of the special form $(x_{\sigma(1)},
  y_{\tau(1)}, x_{\sigma(2)}, y_{\tau(2)}, \ldots, x_{\sigma(n/2)},
  y_{\tau(n/2)})$, for some permutations $\sigma, \tau \in S_{n/2}$.
\end{theorem}

As mentioned in the introduction, this theorem was proven by
connecting loose Hamilton cycles in random 3-uniform hypergraphs with
rainbow Hamilton cycles in randomly edge-colored random graphs, and
applying a result of Janson and Wormald \cite{JW}.  We will use the
reverse connection to demonstrate the equivalence between Theorem
\ref{th1}(a) and Theorem \ref{th:loose-3}.

\vspace{3mm}

\noindent \textbf{Proof of Theorem \ref{th1}(a).}\, Let $K$ be the
constant in Theorem \ref{th:loose-3}.  We are given an even integer
$n$, and a graph $G \sim G_{n,p}$ with $p > \frac{K \log n}{n}$ on
vertex set $\{v_1, \ldots, v_n\}$, whose edges are randomly colored
from the set $\{c_1, \ldots, c_n\}$.  Construct an auxiliary 3-uniform
hypergraph $H$ with vertex set $\{v_1, \ldots v_n, c_1, \ldots,
c_n\}$, by taking the hyperedge $\{v_i, v_j, c_k\}$ whenever the edge
$v_i v_j$ appears in $G$, with color $c_k$.  Note that every such
hyperedge appears independently with probability $\frac{p}{n} >
\frac{K \log n}{n^2}$, since $v_i v_j$ appears in $G$ with probability
$p$, and receives color $c_k$ with probability $\frac{1}{n}$.
Therefore, Theorem \ref{th:loose-3} implies that $H$ has a loose
Hamilton cycle \whp, of the form $(v_{\sigma(1)}, c_{\tau(1)},
v_{\sigma(2)}, c_{\tau(2)}, \ldots, c_{\tau(n/2)})$.  This corresponds
to the Hamilton cycle $(v_{\sigma(1)}, \ldots, v_{\sigma(n/2)})$ in
$G$, with edges colored $c_{\tau(1)}, \ldots, c_{\tau(n/2)}$, hence
rainbow.  \proofend

\section{Proof of Theorem \ref{th1}(b): high level description}
\label{sec:overview}

Let $\e, \th > \frac{100}{\sqrt{\log \log n}}$ be given.  We will
implicitly assume throughout (when convenient) that they are
sufficiently small.  Our proof proceeds in three phases, so our
parameters come in threes.  Let us arbitrarily partition the $\kappa =
(1 + \th) n$ colors into three disjoint groups $C_1 \cup C_2 \cup
C_3$, with sizes
\begin{displaymath}
  |C_1| = \th_1 n,
  \qquad
  |C_2| = (1 + \th_2) n,
  \qquad
  |C_3| = \th_3 n.
\end{displaymath}
We will analyze the random edge generation in three stages, so we
define the probabilities
\begin{displaymath}
  p_1 = \frac{\e_1 \log n}{2 n},
  \qquad
  p_2 = \frac{(1 + \e_2) \log n}{2n},
  \qquad
  p_3 = \frac{\e_3 \log n}{2 n}.
\end{displaymath}
The $\e$'s and $\th$'s are defined by the relations
\begin{equation}
  \label{eq:th-e}
  \e_1 = \e_2 = \e_3 = \frac{\e}{3},
  \qquad
  \th_1 = \th_3 = \min\left\{ \frac{\th}{3}, \frac{\e_2}{4} \right\},
  \qquad
  \th_2 = \th - \th_1 - \th_3.
\end{equation}
(We would have taken $\th_1 = \th_2 = \th_3 = \frac{\th}{3}$, except
that Lemma \ref{lem:cover-S00} requires $\th_1 + \th_3 \leq
\frac{\e_2}{2}$.)

\subsection{Underlying digraph model}
\label{sec:underlying-digraph-model}

It is more convenient for our entire argument to work with directed
graphs, as this will allow us to conserve independence.  Recall that
$\Dnp$ is the model where each of the $n(n-1)$ possible directed edges
appears independently with probability $p$.  We generate a random
colored undirected graph via the following procedure.  First, we
independently generate three digraphs $D_1^\circ = D_{n,p_1}$,
$D_2^\circ = D_{n,p_2}$, and $D_3^\circ = D_{n,p_3}$, and color all of
the directed edges from the full set of colors.

We next use the $D_i^\circ$ to construct a colored undirected graph
$G$, by taking the undirected edge $uv$ if and only if at least one of
$\rt{uv}$ or $\rt{vu}$ appear among the $D_i^\circ$.  The colors of
the undirected edges are inherited from the colors of the directed
edges, in the priority order $D_1^\circ$, $D_2^\circ$, $D_3^\circ$.
Specifically, if $\rt{uv}$ or $\rt{vu}$ appear already in $D_1^\circ$,
then $uv$ takes the color used in $D_1^\circ$ even if $\rt{uv}$ or
$\rt{vu}$ appear again in $D_3^\circ$, say.  In the event that both
$\rt{uv}$ and $\rt{vu}$ appear in $D_1^\circ$, the color of $\rt{uv}$
is used for $uv$ with probability $1/2$, and the color of $\rt{vu}$ is
used otherwise.  Similarly, if neither of $\rt{uv}$ nor $\rt{vu}$
appear in $D_1^\circ$, but $\rt{uv}$ appears in both $D_2^\circ$ and
$D_3^\circ$, the color used in $D_2^\circ$ takes precedence.  It is
clear that the resulting colored graph $G$ has the same distribution
as $G_{n,p,\kappa}$, with
\begin{displaymath}
  p = 1 - (1-p_1)^2 (1-p_2)^2 (1-p_3)^2
  =
  (1 + \e + O(\e^2)) \frac{\log n}{n}.
\end{displaymath}

\subsection{Partitioning by color}
\label{sec:d-to-g}

In each of our three phases, we will use one group of edges and one
group of colors.  Since each $D_i^\circ$ contains edges colored from
the entire set $C_1 \cup C_2 \cup C_3$, for each $i$ we define $D_i
\subset D_i^\circ$ to be the spanning subgraph consisting of all
directed edges whose color is in $C_i$.

Our final undirected graph is generated by superimposing directed
graphs and disregarding the directions.  Consequently, we do not need
to honor the directions when building Hamilton cycles.  To account for
this, we define three corresponding colored undirected graphs $G_1$,
$G_2$, and $G_3$.  These will be edge-disjoint, respecting priority.

The first, $G_1$, is constructed as follows.  For each pair of
vertices $u,v$ with $\rt{uv} \in D_1$ but $\rt{vu} \not \in
D_1^\circ$, place $uv$ in $G_1$ in the same color as $\rt{uv}$.  If
both $\rt{uv}$ and $\rt{vu}$ are in $D_1$, we still place the edge
$uv$ in $G_1$, but randomly select either the color of $\rt{uv}$ or of
$\rt{vu}$.  However, if $\rt{uv} \in D_1$ but $\rt{vu} \in D_1^\circ
\setminus D_1$, then $uv$ is only placed in $G_1$ with probability
$1/2$; if it is placed, it inherits the color of $\rt{uv}$.  Note that
this construction precisely captures all undirected edges arising from
$D_1^\circ$, using colors in $C_1$.

We are less careful with $G_2$, as our argument can afford to discard
all edges that arise from multiply covered pairs.  Specifically, we
place $uv \in G_2$ if and only if $\rt{uv} \in D_2 \setminus
D_1^\circ$ and $\rt{vu} \not\in D_1^\circ \cup D_2^\circ$.  As the
pair $\{u,v\}$ is now spanned by only one directed edge in
$D_2^\circ$, the undirected edge $uv$ inherits that unique color.  We
define $G_3$ similarly, placing $uv \in G_3$ if and only if $\rt{uv}
\in D_3 \setminus (D_1^\circ \cup D_2^\circ)$ and $\rt{vu} \not\in
D_1^\circ \cup D_2^\circ \cup D_3^\circ$.

In this way, we create three edge-disjoint graphs $G_i$.  By our
observations in the previous section, we may now focus on finding an
(undirected) rainbow Hamilton cycle in $G_1 \cup G_2 \cup G_3$.
Importantly, note that in terms of generating colored undirected
edges, the digraph $D_1^\circ$ has higher ``priority'' than
$D_2^\circ$ or $D_3^\circ$.  So, for example, the generation of $G_1$
is not affected by the presence or absence of edges from $D_2^\circ$
or $D_3^\circ$.

\subsection{Main steps}
\label{sec:main-steps}

We generally prefer to work with $G_i$ and $D_i$ instead of
$D_i^\circ$ because we are guaranteed that the edge colors lie in the
corresponding $C_i$.  This allows us to build rainbow segments in
separate stages, without worrying that we use the same color twice.
Let $d_i^+(v)$ denote the out-degree of $v$ in $D_i$.  We now define a
set $S$ of vertices that need special treatment.  We first let
\begin{displaymath}
  S_0 = S_{0,1} \cup S_{0,2} \cup S_{0,3},
\end{displaymath}
where
\begin{align}
& S_{0,1}=\set{v :d_1^+(v)\leq \frac{\e_1\th_1}{20}\log n}\label{S1}\\
& S_{0,2}=\set{v :d_2^+(v)\leq \frac{1}{20}\log n} \label{S2}\\
& S_{0,3}=\set{v :d_3^+(v)\leq \frac{\e_3\th_3}{20}\log n}. \label{S3}
\end{align}
Also, define $\g = \min\set{ \frac{1}{4}, \frac{1}{4} \e_1 \th_1,
  \frac{1}{4} \e_3 \th_3 }$, and note that the constraints on $\e,
\th$ in Theorem \ref{th1} imply the bound
\begin{equation}
  \label{ineq:gamma-lower}
  \g
  >
  \frac{1}{4} \cdot \frac{1}{3} \cdot \frac{1}{12}
  \cdot \bfrac{100}{\sqrt{\log \log n}}^2
  >
  \frac{1}{\log \log n}.
\end{equation}

\begin{lemma}
  \label{lem1}
  The set $S_0$ satisfies $|S_0| \leq \frac{1}{3} n^{1-\g}$ \whp.
\end{lemma}

The vertices in $S_0$ are delicate because they have low degree.  We
also need to deal with vertices having several neighbors in $S_0$.
For this, we define a sequence of sets $S_0,S_1,\ldots,S_t$ in the
following way.  Having chosen $S_t$, if there is still a vertex $v
\not \in S_t$ with at least 4 out-neighbors in $S_t$ (in any of the
graphs $D_1$, $D_2$, or $D_3$), we let $S_{t+1} = S_t \cup \{v\}$ and
continue.  Otherwise we stop at some value $t=T$ and take $S=S_T$.

\begin{lemma}
  \label{lem2}
  With probability $1 - o(n^{-1})$, the set $S$ contains at most
  $n^{1-\g}$ vertices.
\end{lemma}

To take care of the dangerous vertices in $S$, we find a collection of
vertex disjoint paths $Q_1,Q_2,\ldots,Q_s,\,s=|S|$ such that (i) each
path uses undirected edges in $G_2$, (ii) all colors which appear on
these edges are distinct, (iii) all interior vertices of the paths are
vertices of $S$, (iv) every vertex of $S$ appears in this way, and (v)
the endpoints of the paths are not in $S$.  Let us say that these
paths \emph{cover} $S$.

\begin{lemma}
  \label{lem3}
  The graph $G_2$ contains a collection $Q_1,Q_2,\ldots,Q_s$ of paths
  that cover $S$ \whp.
\end{lemma}

The next step of our proof uses a random greedy algorithm to find a
rainbow path of length close to $n$, avoiding all of the previously
constructed $Q_i$.

\begin{lemma}
  \label{lem4}
  The graph $G_2$ contains a rainbow path $P$ of length $n' = n -
  \frac{n}{\sqrt[3]{\log n}}$ \whp.  Furthermore, $P$ is entirely
  disjoint from all of the $Q_i$, and all colors used in $P$ and the
  $Q_i$ are distinct and from $C_2$.
\end{lemma}

Let $U$ be the vertices outside $P$.  (Note that $U$ contains all of
the paths $Q_i$.)  In order to link the vertices of $U$ into $P$, we
split $P$ into short segments, and use the edges of $G_3$ to splice
$U$ into the system of segments.  We will later use the edges of $G_1$
to link the segments back together into a rainbow Hamilton cycle, so
care must be taken to conserve independence.  The following lemma
merges all vertices of $U$ into the collection of segments, and
prepares us for the final stage of the proof.  Here, $d_1^+(v; A)$
denotes the number of $D_1$-edges from $v$ to a set $A$.  Let
\begin{equation}
  \label{eq:L}
  L = \max\left\{
  15 \cdot e^{\frac{40}{\e_3 \th_3}}   
  ,
  \frac{7}{\th_1}
  \right\},
\end{equation}
and note that our conditions on $\e, \th$ in Theorem \ref{th1},
together with \eqref{eq:th-e}, imply that $\e_3 \th_3 > \frac{1}{3}
\cdot \frac{1}{12} \cdot \big( \frac{100}{\sqrt{\log \log n}} \big)^2 >
\frac{277}{\log \log n}$, so we have
\begin{equation}
  \label{ineq:L-upper}
  L
  <
  \max\left\{
    15 \cdot e^{\frac{40}{277} \log \log n}
    ,
    7 \cdot \frac{12 \sqrt{\log \log n}}{100}
  \right\}
  <
  \sqrt[6]{\log n}.
\end{equation}

\begin{lemma}
  \label{lem5}
  With probability $1 - o(1)$, the entire vertex set can be
  partitioned into segments $I_1, \ldots, I_r$, with $r \leq
  \frac{n}{L}$, such that the edges which appear in the segments all
  use different colors from $C_2 \cup C_3$.  The segment endpoints are
  further partitioned into $A \cup B$, with each segment having one
  endpoint in $A$ and one in $B$, such that every $a \in A$ has
  $d_1^+(a; B) \geq \frac{\e_1 \th_1}{200 L} \log n$, and every $b \in
  B$ has $d_1^+(b; A) \geq \frac{\e_1 \th_1}{200 L} \log n$.  All of
  the numeric values $d_1^+(a; B)$ and $d_1^+(b; A)$ have already been
  revealed, but the locations of the corresponding edges are still
  independent and uniform over $B$ and $A$, respectively.
\end{lemma}

The final step links together the segments $I_1, \ldots, I_r$ using
distinctly-colored edges from $G_1$.  For this, we create an auxiliary
colored directed graph $\Gamma$, which has one vertex $w_k$ for each
interval $I_k$.  There is a directed edge $\rt{w_j w_k} \in \Gamma$ if
there is an edge $e \in G_1$ between the $B$-endpoint of $I_j$ and the
$A$-endpoint of $I_k$; it inherits the color of $e$.  Since all colors
of edges in $\Gamma$ are from $C_1$, it therefore suffices to find a
rainbow Hamilton directed cycle in $\Gamma$.  We will find this by
connecting $\Gamma$ with a well-studied random directed graph model.

\begin{definition}
  \label{def:d-in-d-out}
  The $d$-in, $d$-out random directed graph model $D_{d\text{-in},
    d\text{-out}}$ is defined as follows.  Each vertex independently
  chooses $d$ out-neighbors and $d$ in-neighbors uniformly at random,
  and all resulting directed edges are placed in the graph.  Due to
  independence, it is possible that a vertex $u$ selects $v$ as an
  out-neighbor, and $v$ also selects $u$ as an in-neighbor.  In that
  case, instead of placing two repeated edges $\rt{uv}$, place only
  one.
\end{definition}

Instead of proving Hamiltonicity from scratch, we apply the following
theorem of Cooper and Frieze.

\begin{theorem}
  (See \cite{CF-2-inout}, Theorem 1.)
  \label{thm:2-in-2-out}
  The random graph $D_{\text{2-in},\text{2-out}}$ contains a directed
  Hamilton cycle \whp.
\end{theorem}

This result does not take colors into account, however.  Fortunately,
in equation \eqref{eq:L}, we define $L$ to be large enough to allow us
to select a subset of $G_1$-edges which is itself already rainbow.
The analysis of this procedure is the heart of the proof of the final
step.

\begin{lemma}
  \label{lem7}
  The colored directed graph $\Gamma$ contains a rainbow directed
  Hamilton cycle \whp.
\end{lemma}

Since each directed edge of $\Gamma$ corresponds to an undirected
$G_1$-edge from a $B$-endpoint of an interval to an $A$-endpoint of
another interval, a directed Hamilton cycle in $\Gamma$ corresponds to
a Hamilton cycle linking all of the intervals together.  Lemma
\ref{lem7} establishes that it is possible to choose these linking
edges as a rainbow set from $C_1$.  The edges within the intervals
were themselves colored from $C_2 \cup C_3$, so the result is indeed a
rainbow Hamilton cycle in the original graph, as desired.

\section{Proofs of intermediate lemmas}
\label{lemmas}

In the remainder of this paper, we prove the lemmas stated in the
previous section.  Although the first lemma is fairly standard, we
provide all details, and use the opportunity to formally state several
other well-known results which we apply again later.

\subsection{Proof of Lemma \ref{lem1}}

Our first lemma controls the number of vertices whose degrees in $D_i$
are too small.  Recall from Section \ref{sec:underlying-digraph-model}
that the $D_i^\circ$ are independently generated.  Their edges are
then independently colored, and the edges of $D_i^\circ$ which receive
colors from $C_i$ are collected into $D_i$.  (Priorities only take
effect when we form the $G_i$ in Section \ref{sec:d-to-g}.)
Therefore, the out-degrees $d_i^+(v)$ of vertices $v$ in $D_i$ are
distributed as
\begin{align*}
  d_1^+(v) & \sim \bin{n-1,
    p_1 \cdot \frac{\th_1}{1 + \th_1 + \th_2 + \th_3}}
  \geq \bin{0.99 n, \frac{0.49 \e_1 \th_1 \log n}{n}} \\
  d_2^+(v) & \sim \bin{n-1,
    p_2 \cdot \frac{1 + \th_2}{1 + \th_1 + \th_2 + \th_3}}
  \geq \bin{0.99n, \frac{0.49 \log n}{n}} \\
  d_3^+(v) & \sim \bin{n-1,
    p_3 \cdot \frac{\th_3}{1 + \th_1 + \th_2 + \th_3}}
  \geq \bin{0.99n, \frac{0.49 \e_3 \th_3 \log n}{n}}. \\
\end{align*}
Thus the expected size of $S_0$ satisfies
\begin{displaymath}
  \E{|S_0|} \le n(\rho_1 + \rho_2 + \rho_3),
\end{displaymath}
where
\begin{align*}
  \rho_1 &=
  \pr{\bin{0.99 n, \frac{0.49 \e_1 \th_1 \log n}{n}} \leq \frac{\e_1 \th_1 \log n}{20}} \\
  \rho_2 &=
  \pr{\bin{0.99 n, \frac{0.49 \log n}{n}} \leq \frac{\log n}{20}} \\
  \rho_3 &=
  \pr{\bin{0.99 n, \frac{0.49 \e_3 \th_3 \log n}{n}} \leq \frac{\e_3 \th_3 \log n}{20}}.
\end{align*}

We will repeatedly use the following case of the Chernoff lower tail
bound, which we prove with an appropriate explicit constant.

\begin{lemma}
  \label{lem:bin9}
  The following holds for all sufficiently large $mq$, where $m$ is a
  positive integer and $0 < q < 1$ is a real number.
  \begin{displaymath}
    \pr{ \bin{m, q} \leq \frac{1}{9} mq } < e^{-0.533 mq}.
  \end{displaymath}
\end{lemma}

\proofstart Calculation yields
\begin{displaymath}
  \pr{ \bin{m, q} \leq \frac{1}{9} mq }
  =
  \sum_{k=0}^{mq/9} \binom{m}{k} q^k (1-q)^{m-k} \\
  <
  \sum_{k=1}^{mq/9} \bfrac{emq}{k}^k e^{-\frac{8}{9} mq}.
\end{displaymath}
The function $\big( \frac{C}{k} \big)^k = \ex{k(\log C - \log k)}$ is
increasing in $k$ in the range $0 < k < e C$.  Thus
\begin{align*}
  \pr{ \bin{m, q} \leq \frac{1}{9} mq }
  &<
  \frac{mq}{9} \cdot \bfrac{emq}{mq/9}^{mq/9} e^{-\frac{8}{9} mq} \\
  &=
  \frac{mq}{9} \cdot (9e)^{mq/9} e^{-\frac{8}{9} mq} \\
  &=
  e^{ mq ( \frac{1}{9} \log 9e - \frac{8}{9} + o(1)) } \\
  &<
  e^{ -0.533 mq },
\end{align*}
as claimed.
\proofend

\vspace{3mm}

Returning to the proof of Lemma \ref{lem1}, we observe that since
$\frac{1}{20} < \frac{1}{9} \cdot 0.99 \cdot 0.49$, a direct
application of Lemma \ref{lem:bin9} now gives
\begin{align*}
  \rho_2
  &<
  \pr{ \bin{0.99 n, \frac{0.49 \log n}{n}}
    \leq \frac{1}{9} \cdot 0.99 \cdot 0.49 \log n } \\
  &<
  e^{ -0.533 \cdot 0.99 \cdot 0.49 \log n } \\
  &<
  n^{-0.258}.
\end{align*}
A similar argument establishes that $\rho_1 < n^{-0.258 \e_1 \th_1}$
and $\rho_3 < n^{-0.258 \e_3 \th_3}$.  This proves that $\E{|S_0|} =
o(n^{1-\g})$, where we recall our definition $\g = \min\set{
  \frac{1}{4}, \frac{1}{4} \e_2 \th_2, \frac{1}{4} \e_3 \th_3 }$.  We
complete the proof of the lemma by showing that $|S_0|$ is
concentrated around its mean.  For this, we use the Hoeffding-Azuma
martingale tail inequality applied to the vertex exposure martingale
(see, e.g., \cite{AS}).  Recall that a martingale is a sequence $X_0,
X_1, \ldots$ of random variables such that each conditional
expectation $\E{X_{t+1} \mid X_0, \ldots, X_t}$ is precisely $X_t$.

\begin{theorem}
  \label{thm:azuma}
  Let $X_0, \ldots, X_n$ be a martingale, with bounded differences
  $|X_{i+1} - X_i| \leq C$.  Then for any $\lambda \geq 0$,
  \begin{displaymath}
    \pr{X_n \geq X_0 + \lambda}
    \ \leq \
    \exp\left\{-\frac{\lambda^2}{2 C^2 n}\right\}.
  \end{displaymath}
\end{theorem}

Here we consider $|S_0|$ to be a function of $Y_1, Y_2, \ldots, Y_n$
where $Y_k$ denotes the set of edges $\rt{jk}, \rt{kj} \in D_1^\circ
\cup D_2^\circ \cup D_3^\circ$, $j < k$.  The sequence $X_t = \E{|S_0|
  \mid Y_1, \ldots, Y_t}$ is called the vertex-exposure martingale.
There is a slight problem in that the worst-case Lipschitz value for
changing a single $Y_k$ can be too large, while the average case is
good. There are various ways of dealing with this. We will make a
small change in $D^\circ = D_1^\circ \cup D_2^\circ \cup D_3^\circ$.
Let $\hat{D}^\circ$ be obtained from $D^\circ$ by reducing every
degree below $5 \log n$.  We do this in vertex order $v = 1, 2,
\ldots, n$ and delete edges incident with $v$ in descending numerical
order.  We can show that this usually has no effect on $D^\circ$.

\begin{lemma}
  \label{lem:max-degree}
  With probability $1 - o(n^{-1})$, every vertex in $\Gnp$ with $p <
  \frac{1.1 \log n}{n}$ has degree at most $5 \log n$.
\end{lemma}

\proofstart The probability that a single vertex has degree at least
$5 \log n$ is
\begin{align*}
  \pr{\bin{n-1, p} \geq 5 \log n}
  &\leq
  \binom{n}{5\log n} \bfrac{1.1 \log n}{n}^{5 \log n} \\
  &\leq
  \left(
    \frac{e n}{5 \log n}
    \cdot
    \frac{1.1 \log n}{n}
  \right)^{5 \log n} \\
  &=
  \bfrac{1.1e}{5}^{5 \log n} \\
  &= n^{-2.57},
\end{align*}
so a union bound over all vertices gives the result.  \proofend

\vspace{3mm}

Therefore, $\pr{\hat{D}^\circ = D^\circ} = 1 - o(n^{-1})$, and so if
we let $\hat{Z} = |\hat{S}_0|$ be the size of the corresponding set
evaluated in $\hat{D}^\circ$, we obtain $\E{\hat{Z}} = \E{|S_0|} +
o(1) = o(n^{1 - \g})$.  Furthermore, changing a $Y_k$ can only change
$\hat{Z}$ by at most $15 \log n$.  So, we have
\begin{displaymath}
  \pr{\hat{Z}_i \geq \E{\hat{Z}_i} + \frac{1}{4} n^{1-\g}}
  \leq
  \ex{-\frac{n^{2-2\g}/16}{2 (15 \log n)^2 n}}
  <
  o(n^{-1}),
\end{displaymath}
completing the proof of Lemma \ref{lem1}.
\proofend

\subsection{Proof of Lemma \ref{lem2}}
\label{plem2}

We use the following standard estimate to control the densities of
small sets.

\begin{lemma}
  \label{lem:small-sparse}
  With probability $1-o(n^{-1})$, in $\Dnp$ with $p < \frac{\log
    n}{n}$, every set $S$ of fewer than $\frac{4}{e^4} \cdot
  \frac{n}{\log^2 n}$ vertices satisfies $e(S) < 2|S|$.  Here, $e(S)$
  is the number of directed edges spanned by $S$.
\end{lemma}

\proofstart Fix a positive integer $s < \frac{4}{e^4} \cdot
\frac{n}{\log^2 n}$, and consider sets of size $s$.  We may assume
that $s \geq 2$, because a single vertex cannot induce any edges.  The
expected number of sets $S$ with $|S| = s$ and $e(S) \geq 2s$ is at
most
\begin{align*}
  \binom{n}{s} \cdot \binom{s^2}{2s} \bfrac{\log n}{n}^{2s}
  &\leq
  \bfrac{en}{s}^s \cdot \bfrac{es^2}{2s}^{2s} \bfrac{\log n}{n}^{2s} \\
  &=
  \left( s \cdot \frac{e^3 \log^2 n}{4n} \right)^s.
\end{align*}
It remains to show that when this bound is summed over all $2 \leq s <
\frac{4}{e^4} \cdot \frac{n}{\log^2 n}$, the result is still
$o(n^{-1})$.  Indeed, for each $2 \leq s \leq 2 \log n$, the bound is
at most $O \big( \frac{\log^6 n}{n^2} \big)$, so the total
contribution from that part is only $O \big( \frac{\log^7 n}{n^2}
\big) = o(n^{-1})$.  On the other hand, for each $2 \log n < s <
\frac{4}{e^4} \cdot \frac{n}{\log^2 n}$, the bound is at most
\begin{displaymath}
  \left(
    \frac{4}{e^4} \cdot \frac{n}{\log^2 n}
    \cdot
    \frac{e^3 \log^2 n}{4n}
  \right)^{2 \log n}
  =
  \bfrac{1}{e}^{2 \log n}
  =
  n^{-2}.
\end{displaymath}
Thus the total contribution from $2 \log n < s < \frac{4}{e^4} \cdot
\frac{n}{\log^2 n}$ is at most $o(n^{-1})$, as desired.
\proofend

\vspace{3mm}

We are now ready to bound the size of the set $S$ which was created by
repeatedly absorbing vertices with many neighbors in $S_0$.

\vspace{3mm}

\noindent \textbf{Proof of Lemma \ref{lem2}.}\, We actually prove a
stronger statement, which we will need for Lemma \ref{lem:few-nbr-S}.
Suppose we have an initial $S_0$ satisfying $|S_0| < \frac{1}{3}
n^{1-\g}$, as ensured by Lemma \ref{lem1}.  Consider a sequence $S_0',
S_1', S_2', \ldots$ where $S_{t+1}'$ is obtained from $S_t'$ by adding
a vertex $v \notin S_t'$ for which $d_i^+(v; S_t') \geq 3$ for some
$i$.  Note that when this process stops, the final set $S'$ will
contain the set $S$ which our definition obtained by adding vertices
with degree at least 4 into previous $S_t$.

So, suppose for contradiction that this process continues for so long
that some $|S_t'|$ reaches $n^{1-\g} = o\big( \frac{n}{\log^2 n}
\big)$.  Note that $t \geq \frac{2}{3} n^{1-\g}$.  Since each step
introduces at least 3 edges, we must have $e(S_t) \geq 3t \geq 2
n^{1-\g} = 2 |S_t|$.  Yet by construction, $D_1 \cup D_2 \cup D_3$ is
an instance of $D_{n,q}$ for some $q < \frac{\log n}{n}$, so this
contradicts Lemma \ref{lem:small-sparse}.  \proofend

\subsection{Proof of Lemma \ref{lem3}}

In this section, we show that for each vertex $v \in S$, we can find a
disjoint $G_2$-path $Q \ni v$ which starts and ends outside $S$.  We
also need all colors appearing on these edges to be different.  Since
we are working in a regime where degrees can be very small, we need to
accommodate the most delicate vertices first.  Specifically, let
$S_{0,0}$ be the set of all vertices with $d_2(v) \leq \frac{1}{10}
\log n$, where $d_2(v)$ is the degree of $v$ in $G_2$.  Although
$S_{0,0}$ will typically not be entirely contained within $S_{0,2}$,
we can show that it is still usually quite small.

\begin{lemma}
  \label{lem:S00}
  We have $|S_{0,0}| < n^{0.48}$ \whp.
\end{lemma}

\proofstart By construction, $G_2 \sim G_{n,q_2}$, where
\begin{equation}
  q_2
  =
  2 p_2 (1-p_2)
  \cdot
  \frac{1 + \th_2}{1 + \th_1 + \th_2 + \th_3}
  \cdot
  (1 - p_1)^2,
  \label{eq:q2}
\end{equation}
because the first factor is the probability that exactly one of
$\rt{uv}$ or $\rt{vu}$ appears in $D_2^\circ$, the second factor is
the probability that it receives a color from $C_2$, and the third
factor is the probability that neither $\rt{uv}$ nor $\rt{vu}$ appear
in $D_1^\circ$.  Hence for a fixed vertex $v$, its relevant degree in
$G_2$ is distributed as
\begin{displaymath}
  d_2(v) \sim \bin{n-1, q_2}
  \geq \bin{0.99 n, \frac{0.99 \log n}{n}}.
\end{displaymath}
Since $\frac{1}{10} < \frac{1}{9} \cdot 0.99 \cdot 0.99$, Lemma
\ref{lem:bin9} implies that
\begin{align}
  \nonumber
  \pr{ d_2(v) \leq \frac{1}{10} \log n }
  &<
  \pr{ \bin{0.99 n, \frac{0.99 \log n}{n}}
    \leq \frac{1}{9} \cdot 0.99 \cdot 0.99 \log n } \\
  \nonumber &<
  e^{ -0.533 \cdot 0.99 \cdot 0.99 \log n } \\
  &<
  n^{-0.522}. \label{d1real-low}
\end{align}
Therefore, $\E{|S_{0,0}|} < n \cdot n^{-0.522}$, and Markov's
inequality yields the desired result.
\proofend

\vspace{3mm}

We have shown that vertices of $S_{0,0}$ are few in number.  Our next
result shows that they are also scattered far apart.  This will help
us when we construct the covering paths, by preventing paths from
colliding.

\begin{lemma}
  \label{lem:dist5}
  Let $\dist_2(v,w)$ denote the distance between $v$ and $w$ in $G_2$.
  Then, \whp, every pair $v,w \in S_{0,0}$ satisfies $\dist_2(v,w)
  \geq 5$.
\end{lemma}

\proofstart Recall that $G_2 \sim G_{n,q_2}$ with $q_2$ defined as in
\eqref{eq:q2}.  Consider a fixed pair of vertices $v,w$.  For a fixed
sequence of $k \leq 4$ intermediate vertices $x_1, x_2, \ldots, x_k$,
let us bound the probability $q$ that $v,w$ both have $d_2 \leq
\frac{1}{10} \log n$, and all the edges $v x_1, x_1 x_2, x_2 x_3,
\ldots, x_k w$ appear in $G_2$.  First expose the edges $v x_1$, $x_1
x_2$, \ldots, $x_k w$, and then expose the edges between $v$ and $[n]
\setminus \{v, x_1, \ldots, x_k, w\}$, and between $w$ and that set.

This gives the following bound on our probability $q$:
\begin{equation}
  \label{ineq:dist10}
  q \leq \bfrac{1.01 \log n}{n}^{k+1} \cdot
  \pr{ \bin{n - 2 - k, \frac{0.99 \log n}{n}}
    \leq \frac{\log n}{10} }^2
\end{equation}
A calculation analogous to \eqref{d1real-low} bounds the Binomial
probability by $n^{-0.52}$, so taking a union bound over all $O(n^k)$
choices for the $x_i$, for all $0 \leq k \leq 4$, we find that for
fixed $v,w$, the probability that $\dist_2(v,w) < 5$ is at most
\begin{displaymath}
  \sum_{k=0}^4 O(n^k) \cdot \bfrac{1.01 \log n}{n}^{k+1} \cdot
  \left( n^{-0.52} \right)^2
  <
  n^{-2.04 + o(1)}.
\end{displaymath}
Therefore, a final union bound over the $O(n^2)$ choices for $v,w$
completes the proof.
\proofend

\vspace{3mm}

The previous result will help us cover vertices in $S_{0,0}$ with
$G_2$-paths.  However, the objective of this section is to cover all
vertices of $S$.  Although the analogue of Lemma \ref{lem:dist5} does
not hold for $S$, it is still possible to prove that $S$ is sparsely
connected to the rest of the graph.  Recall from
\eqref{ineq:gamma-lower} that $\g = \min\set{ \frac{1}{4}, \frac{1}{4}
  \e_1 \th_1, \frac{1}{4} \e_3 \th_3 } > \frac{1}{\log \log n}$.

\begin{lemma}
  \label{lem:few-nbr-S}
  With respect to edges of $G_2$, every vertex $v$ is adjacent to at
  most $\frac{2}{\g}$ vertices in $S$ \whp.  (This applies whether or
  not $v$ itself is in $S$.)
\end{lemma}

\proofstart Fix a vertex $v$.  Let $S'$ be the set obtained by
constructing the analogous sequences to $S_{0,1}'$, $S_{0,2}'$,
$S_{0,3}'$, $S_1'$, \ldots\ on the graph induced by $[n] \setminus
\{v\}$, where $S_{0,i}'$ are defined as in \eqref{S1}--\eqref{S3}, but
the $S_{t+1}'$ are obtained by adding vertices with at least 3 (not 4)
$D_i$-out-neighbors in $S_t'$.  Clearly, $S'$ contains $S \setminus
\{v\}$, because the effect of ignoring $v$ is compensated for by using
3 instead of 4.  The advantage of using $S'$ instead of $S$ is that
$S'$ can be generated without exposing any edges incident to $v$.  As
we will take a final union bound over the $n$ choices of $v$, it
therefore suffices to show that with probability $1 - o(n^{-1})$, the
particular vertex $v$ has at most $\frac{2}{\g}$ neighbors in $S'$.

For this, we expose all edges of $D_1^\circ \cup D_2^\circ \cup
D_3^\circ$ that are spanned by $[n] \setminus \{v\}$.  Recall that our
proof of Lemma \ref{lem2} already absorbed vertices with 3
out-neighbors (instead of 4), so we have $|S'| \leq n^{1-\g}$ with
probability $1 - o(n^{-1})$.  It remains to control the number of
edges between $v$ and $S'$, so we now expose all edges of $D_1^\circ
\cup D_2^\circ \cup D_3^\circ$.  The $G_2$-edges there appear
independently with probability $q_2$ as defined in \eqref{eq:q2}, so
the probability that at least $\frac{2}{\g}$ edges appear is at most
\begin{align}
  \nonumber \pr{\bin{n^{1-\g}, q_2} \geq \frac{2}{\g}}
  &\leq
  \binom{ n^{1-\g}}{2/\g }
  \bfrac{1.01 \log n}{n}^{2/\g} \\
  \nonumber &\leq
  \left(
    \frac{e \cdot n^{1-\g}}{2/\g}
    \cdot
    \frac{1.01 \log n}{n}
  \right)^{2/\g} \\
  \nonumber &<
  \left(
    \frac{2 \g \log n}{n^\g}
  \right)^{2/\g} \\
  &<
  \frac{(2 \log n)^{2 \log \log n}}{n^2}
  =
  o(n^{-1}).  \label{ineq:few-nbr-S}
\end{align}
Taking a final union bound over all initial choices for $v$ completes
the proof.  \proofend

\vspace{3mm}

We will cover each vertex $v \in S$ with a $G_2$-path by joining two
$G_2$-paths of length up to to 2, each originating from $v$.  It is
therefore convenient to extend the previous result by one further
iteration.

\begin{corollary}
  \label{cor:few-dist2-S}
  With respect to edges of $G_2$, every vertex $v$ is within distance
  two of at most $\big( \frac{2}{\g} \big)^2$ vertices in $S$ \whp.
  (This applies whether or not $v$ itself is in $S$.)
\end{corollary}

\proofstart Fix a vertex $v$.  Construct $S'$ in the same way as in
the proof of Lemma \ref{lem:few-nbr-S}, exposing only edges spanned by
$[n] \setminus \{v\}$.  Using only those exposed edges, let $T \subset
[n] \setminus \{v\}$ be the set of all vertices in $S'$ or adjacent to
$S'$ via edges from $G_2$.  By Lemma \ref{lem:max-degree}, the maximum
degree of $G_2 \setminus \{v\}$ is at most $5 \log n$ with probability
$1 - o(n^{-1})$, so $|T| \leq n^{1-\g} \cdot 5 \log n$.  A similar
calculation to \eqref{ineq:few-nbr-S} then shows that with probability
$1 - o(n^{-1})$, $v$ has at most $\frac{2}{\g}$ neighbors in $T$.

Taking a union bound over all $v$, and combining this with Lemma
\ref{lem:few-nbr-S}, we conclude that \whp, every vertex has at most
$\frac{2}{\g}$ neighbors in $S \cup N(S)$, and each of them has at
most $\frac{2}{\g}$ neighbors in $S$.  This implies the result.
\proofend

\vspace{3mm}

We are now ready to start covering the vertices of $S$ with disjoint
rainbow $G_2$-paths.  The most delicate vertices are those in
$S_{0,0}$, because by definition all other vertices already have
$G_2$-degree at least $\frac{1}{10} \log n$.  Naturally, we take care
of $S_{0,0}$ first.

\begin{lemma}
  \label{lem:cover-S00}
  The colored graph $G_2$ contains a rainbow collection $Q_1, Q_2,
  \ldots, Q_{s'}$ of disjoint paths that cover $S \cap S_{0,0}$ \whp.
\end{lemma}

\proofstart We condition on the high-probability events in Lemmas
\ref{lem:S00}, \ref{lem:dist5}, and \ref{lem:few-nbr-S}, and use a
greedy algorithm to cover each $v \in S \cap S_{0,0}$ with a path of
length 2, 3, or 4.  Recall that $G_2 \sim G_{n,q_2}$, where $q_2$ was
specified in \eqref{eq:q2}.  We can bound $q_2$ by
\begin{displaymath}
  q_2
  >
  (1 + \e_2) \frac{\log n}{n}
  \left( 1 - \frac{\log n}{n} \right)
  \cdot
  (1 - \th_1 - \th_3).
  \left( 1 - \frac{\log n}{n} \right)^2.
\end{displaymath}
Since equation \eqref{eq:th-e} ensures that $\th_1 + \th_3 \leq
\frac{\e_2}{2}$, and our conditions on Theorem \ref{th1} force $\e_2 >
\frac{1}{3} \cdot \frac{100}{\sqrt{\log \log n}} \gg \frac{\log \log
  n}{\log n}$, the minimum degree in $G_2$ is at least two \whp.
Condition on this as well.

Now consider a vertex $v \in S \cap S_{0,0}$, and let $x_1$ and $x_2$
be two of its neighbors.  If both $x_i$ are already outside $S$, then
we use $x_1 v x_2$ to cover $v$.  Otherwise, suppose that $x_1$ is
still in $S$.  Since we conditioned on vertices in $S_{0,0}$ being
separated by distances of at least 5 (Lemma \ref{lem:dist5}), $x_1$
cannot be in $S_{0,0}$, so it has at least $\frac{1}{10} \log n$
$G_2$-neighbors.  These cannot all be in $S$, because we conditioned
on the fact that every vertex has fewer than $\frac{2}{\g} < 2 \log
\log n$ neighbors in $S$ (Lemma \ref{lem:few-nbr-S}).  So, we can pick
one, say $y_1$, such that $y_1 x_1 v$ is a path from outside $S$ to
$v$.  A similar argument allows us to continue the path from $v$ to a
vertex outside $S$ in at most two steps.  Therefore, there is a
collection of paths of length 2--4 covering each vertex in $S \cap
S_{0,0}$.  They are all disjoint, since we conditioned on vertices of
$S_{0,0}$ being separated by distances of at least 5.

At this point, we have exposed all $G_2$-edges spanned by $S$ and its
neighbors, but the only thing we have revealed about their colors is
that they are all in $C_2$.  Now expose the precise colors on all
edges of these paths.  Since we conditioned on $|S_{0,0}| < n^{0.48}$
(Lemma \ref{lem:S00}), the total number of edges involved is at most
$4 \cdot n^{0.48} < n^{0.49}$.  The number of colors in $C_2$ is $(1 +
\th_2) n$, so by a simple union bound the probability that some pair
of edges receives the same color in $C_2$ is at most
\begin{displaymath}
  \binom{n^{0.49}}{2} \frac{1}{(1 + \th_2) n} = o(1).
\end{displaymath}
Therefore, the covering paths form a rainbow set \whp, as desired.
\proofend

\vspace{3mm}

We have now covered the most dangerous vertices of $S$.  The remainder
of this section provides our argument which covers all other vertices
in $S$.

\vspace{3mm}

\noindent \textbf{Proof of Lemma \ref{lem3}.}\, Condition on the
high-probability events of Lemmas \ref{lem2}, \ref{lem:few-nbr-S},
\ref{lem:cover-S00}, and Corollary \ref{cor:few-dist2-S}.  We have
already covered all vertices in $S \cap S_{0,0}$ with disjoint rainbow
paths of lengths up to four (Lemma \ref{lem:cover-S00}).  We cover the
rest of the vertices in $S \setminus S_{0,0}$ with paths of length
two, using a simple iterative greedy algorithm.  Indeed, suppose that
we are to cover a given vertex $v \in S \setminus S_{0,0}$.  Since it
is not in $S_{0,0}$, it has $G_2$-degree at least $\frac{1}{10} \log
n$, and at most $\frac{2}{\g}$ of these neighbors can be within $S$
(Lemma \ref{lem:few-nbr-S}).

Furthermore, we can show that at most $2 (\frac{2}{\g})^2$ of $v$'s
neighbors outside $S$ can already have been used by covering paths.
Indeed, for each neighbor $w \not \in S$ of $v$ which was used by a
previous covering path, we could identify a vertex $x \in S$ adjacent
to $w$ which was part of that covering path.  Importantly, $x$ is
within distance two of $v$, so the collection of all $x$ obtainable in
this way is of size at most $(\frac{2}{\g})^2$, as we conditioned on
Corollary \ref{cor:few-dist2-S}.  Since every covering path uses
exactly two vertices outside $S$, the total number of such $w$ is at
most $2 (\frac{2}{\g})^2$.  Putting everything together, we conclude
that the number of usable $G_2$-edges emanating from $v$ is at least
\begin{displaymath}
  \frac{1}{10} \log n - \frac{2}{\g} - 2 \bfrac{2}{\g}^2
  >
  \frac{1}{11} \log n.
\end{displaymath}
Expose the colors (necessarily from $C_2$) which appear on these
$G_2$-edges.  Of the total of $(1 + \th_2) n$ available, we only need
to avoid at most $4 |S|$ which have already been used on previous
covering paths.  Since we conditioned on $|S| \leq n^{1-\g}$ (Lemma
\ref{lem2}), this is at most $4 n^{1-\g}$ colors to avoid.  We only
need to have two new colors to appear among this collection in order
to add a new rainbow path of length two covering $v$.  Taking another
union bound, we find that the probability that at most one new color
appears is at most
\begin{align*}
  (1 + \th_2)n \cdot \bfrac{4 n^{1-\g} + 1}{(1 + \th_2)n}^{\frac{1}{11} \log n}
  &= o(n^{-1}).
\end{align*}
Here, the first factor of $(1+\th_2)n$ corresponds to the number of
ways to choose the new color to add (or none at all).  Since we only
run our algorithm for $o(n)$ iterations (once per vertex in $S
\setminus S_{0,0}$), we conclude that \whp\ we can cover all vertices
of $S$ with disjoint rainbow $G_2$-paths.  \proofend

\subsection{Proof of Lemma \ref{lem4}}
\label{sec:long-path}

In this section, we construct a rainbow $G_2$-path which contains most
of the vertices of the graph, but avoids all covering paths from the
previous section.  In order to carefully track the independence and
exposure of edges, recall from Section \ref{sec:d-to-g} that $G_2$ is
deterministically constructed from the random directed graphs
$D_1^\circ$, $D_2^\circ$, and $D_3^\circ$.  Let us consider the
generation of the $D_i^\circ$ to be as follows.  The probability that
the directed edge $\rt{vw}$ appears in $D_1$ is $p_1 \cdot
\frac{\th_1}{1 + \th_1 + \th_2 + \th_3}$, so we expose each
$D_1$-out-degree $d_1^+(v)$ by independently sampling from the
$\bin{n-1, p_1 \cdot \frac{\th_1}{1 + \th_1 + \th_2 + \th_3}}$
distribution.  Importantly, we do not reveal the locations of the
out-neighbors.  Similarly, for $D_2$ and $D_3$, we expose all
out-degrees $d_2^+(v)$ and $d_3^+(v)$, each sampled from the
appropriate Binomial distribution.  By Lemma \ref{lem:max-degree}, all
$d_i^+(v) \leq 5 \log n$ \whp; we condition on this.

Note that from this information, we can later fully generate (say)
$D_1$ and $D_1^\circ$ as follows.  At each vertex $v$, we
independently choose $d_1^+(v)$ out-neighbors uniformly at random.
This will determine all $D_1$-edges.  Next, for every edge which is
not part of $D_1$, independently sample it to be part of $D_1^\circ
\setminus D_1$ with probability $p_1 \big( 1 - \frac{\th_1}{1 + \th_1
  + \th_2 + \th_3} \big)$.  This will determine all $D_1^\circ$ edges,
and a similar system will determine all edges of $D_2^\circ$ and
$D_3^\circ$.

Returning to the situation where only the $d_i^+(v)$ have been
exposed, we then construct the $S_{0,i}$ by collecting all vertices
whose $d_i^+(v)$ are too small, and build the sequence $S_0, S_1, S_2,
\ldots, S_t$.  In each iteration of that process, we go over all
vertices which are not yet in the current $S_t$.  At each $v$, we
expose all $D_i$-edges incident to $S_t$.  For this section, we will
only care about the $D_2$-out-edges from $v \not\in S$ (initially
counted by $d_2^+(v)$) that are \emph{not}\/ consumed in this process.
Fortunately, at each exposure stage, there is a clear distribution on
the number of these out-edges that are consumed toward $S_t$, and this
will only affect the number, not the location, of the out-edges which
are not consumed.

Therefore, after this procedure terminates, we will have a final set
$S$, and the set of revealed (directed) edges is precisely those edges
spanned by $S$, together with all those between $S$ and $V_1 = [n]
\setminus S$.  This set of revealed edges is exactly what is required
to construct the covering paths $Q_1, \ldots, Q_s$ in Lemma
\ref{lem3}.  Within $V_1$, the precise locations of the edges are not
yet revealed.  Instead, for each vertex $v \in V_1$, there is now a
number $d_2^*(v)$, corresponding to the number of $D_2$-out-edges from
$v$ to vertices outside $S$.

We now make two crucial observations.  First, the distributions of
where these endpoints lie are still independent and uniform over
$V_1$.  Second, every $d_2^*(v) \leq 5 \log n$ and $d_2^*(v) \geq
\frac{1}{20} \log n - 3 > \frac{1}{21} \log n$, because if there were
4 out-edges from $v$ to $S$, then $v$ should have been absorbed into
$S$ during the process.

This abundance of independence makes it easy to analyze a simple
method for finding a long path, based on a greedy algorithm with
backtracking.  (This procedure is similar to that used in \cite{DLV}
by Fernandez de la Vega.)  Indeed, the most straightforward attempt
would be to start building a path, and at each iteration expose the
out-edges of the final endpoint, as well as their colors.  If there is
an option which keeps the path rainbow, we would follow that edge, and
repeat.  If not, then we should backtrack to the latest vertex in the
path which still has an option for extension.

We formalize this in the following algorithm.  Particularly dangerous
vertices will be coded by the color red (not related to the colors of
the edges in the $\Gnpk$).  Let $V_2 \subset V_1$ be the set of all
vertices which are not involved in the covering paths $Q_i$.  We will
find a long $G_2$-path within $V_2$ which avoids all of the covering
paths.

\vspace{3mm}

\noindent \textbf{Algorithm.}
\begin{enumerate}
\item Initially, let all vertices of $V_2$ be uncolored, and select an
  arbitrary vertex $v \in V_2$ to use as the initial path $P_0 =
  \{v\}$.  Let $U_0 = V_2 \setminus \{v\}$.  This is the set of
  ``untouched'' vertices.  Let $R_0 = \emptyset$.  This will count the
  ``red'' vertices.

\item Now suppose we are at time $t$.  If $|U_t| < \frac{n}{2
    \sqrt[3]{\log n}}$, terminate the algorithm.

\item If the final endpoint $v$ of $P_t$ is not red, then expose the
  first $\frac{1}{2} d_2^*(v)$ of $v$'s $D_2$-out-neighbors.  If none of
  them lies in $U_t$, via an edge color not yet used by $P_t$ or any
  of the covering paths $Q_i$, then color $v$ red, setting $U_{t+1} =
  U_t$, $P_{t+1} = P_t$, and $R_{t+1} = R_t \cup \{v\}$.

  Otherwise, arbitrarily choose one of the suitable out-neighbors $w
  \in U_t$.  Set $U_{t+1} = U_t \setminus \{w\}$.  Expose whether
  $\rt{vw} \in D_1^\circ$, $\rt{wv} \in D_1^\circ$, or $\rt{wv} \in
  D_2^\circ$.  If none of those three directed edges are present, then
  add $w$ to the path, setting $P_{t+1} = P_t \cup \{w\}$ and $R_{t+1}
  = R_t$.  Otherwise, color both $v$ and $w$ red, and set $P_{t+1} =
  P_t$ and $R_{t+1} = R_t \cup \{v, w\}$.

\item If the final endpoint $v$ of $P_t$ is red, then expose the
  second $\frac{1}{2} d_2^*(v)$ of $v$'s $D_2$-out-neighbors.  First
  suppose that none of them lies in $U_t$, via an edge color not yet
  used by $P_t$ or any of the covering paths $Q_i$.  In this case,
  find the last vertex $v'$ of $P_t$ which is not red, color it red,
  and make it the new terminus of the path.  That is, set $U_{t+1} =
  U_t$, let $P_{t+1}$ be $P_t$ up to $v'$, and set $R_{t+1} = R_t \cup
  \{v'\}$.  If $v'$ did not exist (i.e., all vertices of $P_t$ were
  already red), then instead let $v'$ be an arbitrary vertex of $U_t$
  and restart the path, setting $P_{t+1} = \{v'\}$, $R_{t+1} = R_t$,
  $U_{t+1} = U_t \setminus \{v'\}$.

  On the other hand, if $v$ has a suitable out-neighbor $w \in U_t$,
  then set $U_{t+1} = U_t \setminus \{w\}$.  Expose whether $\rt{vw}
  \in D_1^\circ$, $\rt{wv} \in D_1^\circ$, or $\rt{wv} \in D_2^\circ$.
  If none of those three directed edges are present, then add $w$ to
  the path, setting $P_{t+1} = P_t \cup \{w\}$ and $R_{t+1} = R_t$.
  Otherwise, color $w$ red, find the last vertex $v'$ of $P_t$ which
  is not red, and follow the remainder of the first paragraph of this
  step.
\end{enumerate}

\vspace{3mm}

\noindent The key observation is that the final path $P_T$ contains
every non-red vertex which lies in $V_2 \setminus U_T$.  Since Lemmas
\ref{lem2} and \ref{lem3} imply that $|V_2| \geq n - 3 |S| \geq n - 3
n^{1-\g}$, and we run until $U_T < \frac{n}{2 \sqrt[3]{\log n}}$,
Lemma \ref{lem4} therefore follows from the following bound.

\begin{lemma}
  \label{lem:algo-red-whp}
  The final number of red vertices is at most $n \cdot
  e^{-\frac{1}{300} \sqrt[3]{\log n}}$ \whp.
\end{lemma}

\proofstart The color red is applied in only two situations.  The
first is when we expose whether any of $\rt{vw} \in D_1^\circ$,
$\rt{wv} \in D_1^\circ$, or $\rt{wv} \in D_2^\circ$ hold.  To expose
whether $\rt{vw} \in D_1^\circ$, we reveal whether $\rt{vw} \in D_1$,
using the previously exposed value of $d_1^*(v)$, which we already
conditioned on being at most $5 \log n$.  Since $v$'s
$D_1$-out-neighbors are uniform, the probability that $\rt{vw} \in
D_1$ is at most $\frac{5 \log n}{(1 - o(1))n}$.  If it is not in
$D_1$, the probability that it is in $D_1^\circ \setminus D_1$ is
bounded by $\frac{\log n}{n}$ by the description at the beginning of
Section \ref{sec:long-path}.  The analysis for the other two cases are
similar, so a union bound gives that the chance that any of $\rt{vw}
\in D_1^\circ$, $\rt{wv} \in D_1^\circ$, or $\rt{wv} \in D_2^\circ$
hold is at most $3 \cdot \frac{7 \log n}{n}$.  Note that this occurs
at most $n$ times, because each instance reduces the size of $U_t$ by
1.  Hence the expected number of red vertices of this type is at most
$O(\log n)$, which is of much smaller order than $n
e^{-\Theta(\sqrt[3]{\log n})}$.

The other situation in which red is applied comes immediately after
the failed exposure of some $k = \frac{1}{2} d_2^*(v) > \frac{1}{42}
\log n$ $D_2$-out-neighbors, in either of Steps 3 or 4.  Failure means
that all $k$ of them either fell outside $U_t$, or had edge colors
already used in $P_t$ or some covering path $Q_i$.  Step 2 controls
$|U_t| \geq \frac{n}{2 \sqrt[3]{\log n}}$, and the total number of
colors used in $P_t$ or any covering path $Q_i$ is at most $n -
|U_t|$, out of the $(1+\th_2) n$ available.  Further note that because
of our order of exposure, there is a set $T$ of size at most $3 \log
n$ such that $v$'s $D_2$-out-neighbors are uniformly distributed over
$V_1 \setminus T$.  This is because we have exposed whether $\rt{vu}$
was a $D_2$-edge, for the predecessor $u$ of $v$ along $P_t$, and we
may also have already exposed the first half of $v$'s
$D_2$-out-neighbors in a prior round, which could consume up to
$\frac{1}{2} \cdot 5 \log n$ vertices.  Therefore, the chance that a
given out-neighbor exposure is successful (i.e., lands inside $U_t$,
via one of the $\geq |U_t|$ unused colors), is at least
\begin{displaymath}
  \frac{|U_t \setminus T|}{|V_1 \setminus T|}
  \cdot
  \frac{|U_t|}{(1 + \th_2)n}
  \geq
  \left(
    \frac{n}{2 \sqrt[3]{\log n}}
    \cdot
    \frac{1}{(1-o(1))n}
  \right)
  \cdot
  \left(
    \frac{n}{2 \sqrt[3]{\log n}}
    \cdot
    \frac{1}{(1 + \th_2)n}
  \right)
  >
  \frac{1}{5 (\log n)^{2/3}}.
\end{displaymath}
We conclude that the chance that all $k \geq \frac{1}{42} \log n$ fail
is at most
\begin{displaymath}
  (1+o(1))\left(
    1 - \frac{1}{5 (\log n)^{2/3}}
  \right)^{\frac{1}{42} \log n}
  <
  e^{ -\frac{1}{210} \sqrt[3]{\log n} }.
\end{displaymath}
Since we will not perform this experiment more than twice for each of
the $n$ vertices, linearity of expectation and Markov's inequality
imply that \whp, the final total number of red vertices is at most
$n \cdot e^{- \frac{1}{300} \sqrt[3]{\log n}}$,
as desired.
\proofend

\subsection{Proof of Lemma \ref{lem5}}

\begin{figure}[b]
  \centering
  \includegraphics[scale=0.8]{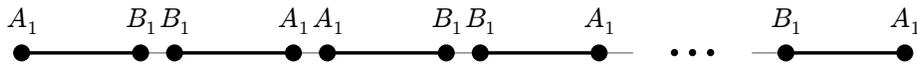}
  \caption{\footnotesize The long path $P$, divided into consecutive
    intervals of length $L$.  Endpoints of successive intervals are
    adjacent via original edges of $P$.  The set of interval endpoints
    has been partitioned into $A_1 \cup B_1$.  Note that endpoints
    which are adjacent via an edge of $P$ are always assigned to the
    same set.}
  \label{fig:path-AB}
\end{figure}

At this point, we have a rainbow $G_2$-path $P$ of length $n' \geq n -
\frac{n}{\sqrt[3]{\log n}}$, which is disjoint from the paths $Q_i$
which cover $S$.  Recall from \eqref{eq:L} and \eqref{ineq:L-upper}
that we defined $L = \max\big\{ 15 e^{40/(\e_3 \th_3)} ,
\frac{7}{\th_1} \big\} < \sqrt[6]{\log n}$.  Split $P$ into $r =
\frac{n'}{L}$ segments of length $L$, as in Figure \ref{fig:path-AB}.
If $n'$ is not divisible by $L$, we may discard the remainder of $P$,
because $L < \sqrt[6]{\log n}$.

Partition the $2 r$ endpoints into two sets $A_1 \cup B_1$ so that
each segment has one endpoint in each set, but there are no vertices
$a \in A_1$ and $b \in B_1$ which are consecutive along $P$.  By
possibly discarding the final interval (which will only cost an
additional $L < \sqrt[6]{\log n}$), we may ensure that the initial and
final endpoints are both in $A_1$.

The reason for our unusual partition is as follows.  In our
construction thus far, we already needed to expose the locations of
some $D_1^\circ$-edges, since they had priority over the $D_2^\circ$.
In certain locations, we have revealed that there are no
$D_1^\circ$-edges.  In particular, between every consecutive pair of
vertices $u, v$ on the path $P$, we found a $D_2$-edge, and confirmed
the absence of any $D_1^\circ$-edges.

Fortunately, our construction did not expose any $D_1^\circ$-edges
between non-consecutive vertices of the path $P$.  In particular, if
we now wished, for any vertex $a \in A_1$, we could expose the number
$N$ of its $D_1$-out-neighbors that lie in $B_1$; then, the
distribution of these $N$ out-neighbors would be uniform over $B_1$.
This uniformity is crucial, and would not hold, for example, if some
vertex of $B_1$ were consecutive with $a$ along $P$.

The proof of Lemma \ref{lem5} breaks into the following steps.  Recall
that $d_1^+(v; T)$ denotes the number of $D_1$-edges from a vertex $v$
to a subset $T$ of vertices.  Say that $v$ is \emph{$T$-good} if
$d_1^+(v; T) \geq \frac{\e_1 \th_1}{180 L} \log n$; call it
\emph{$T$-bad} otherwise.

\begin{figure}[p]
  \centering
  \includegraphics[scale=0.8]{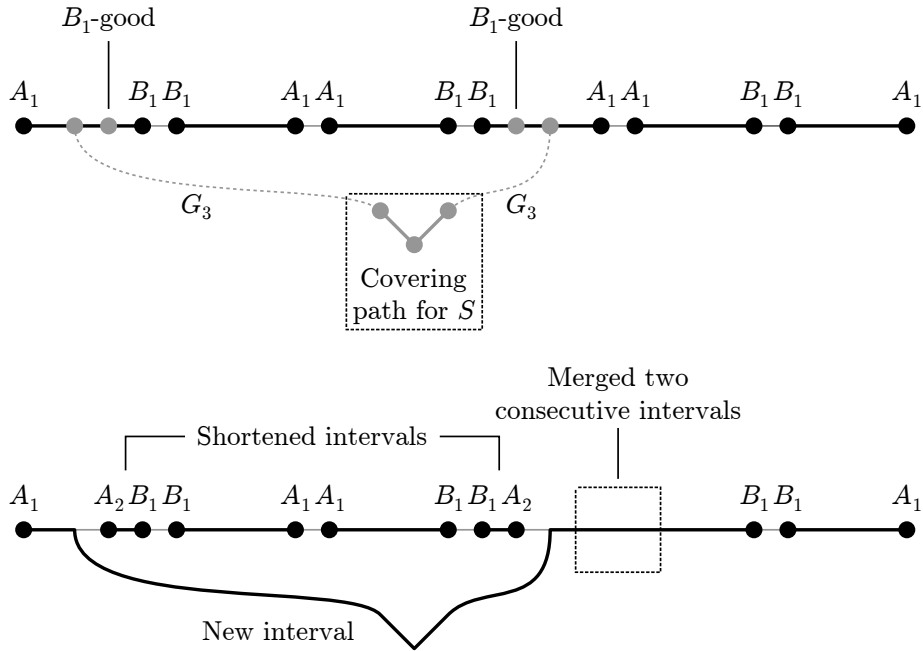}
  \caption{\footnotesize A covering path of $S$ is absorbed into the
    system of intervals using $G_3$-edges.  Note that the resulting
    endpoint partition still has one $A$-endpoint and one $B$-endpoint
    in every interval.  Importantly, the direction of the splicing is
    such that all new endpoints are of type-$A$, as indicated by the
    $A_2$-vertices.  This is why we continue the new interval
    rightward, through to the next $B$-endpoint.}
  \label{fig:path-absorb}
\end{figure}

\begin{figure}[p]
  \centering
  \includegraphics[scale=0.8]{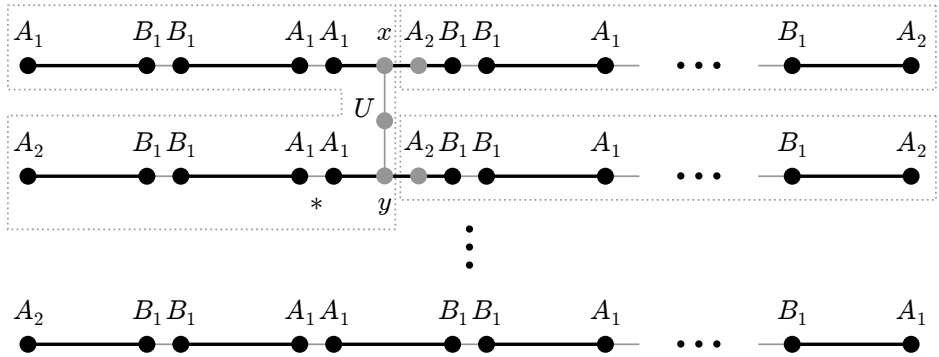}
  \caption{\footnotesize Evolution of block partition during
    absorption.  Each horizontal row represents an original block,
    within which the successive segments have their endpoints
    connected by edges of the original path $P$.  The vertical gray
    path from $x$ to $y$ represents the absorption of a new vertex
    into the collection of segments, involving two different blocks.
    This operation cuts the two edges between $x, y$ and their
    adjacent $A_2$-vertices, and adds back the $P$-edge marked by the
    asterisk.  Afterward, the segments can be re-partitioned into new
    blocks (see the gray dotted lines), with all initial and final
    endpoints in each block of type-$A$, and $B_1$-good. }
  \label{fig:path-strings}
\end{figure}

\begin{figure}[p]
  \centering
  \includegraphics[scale=0.8]{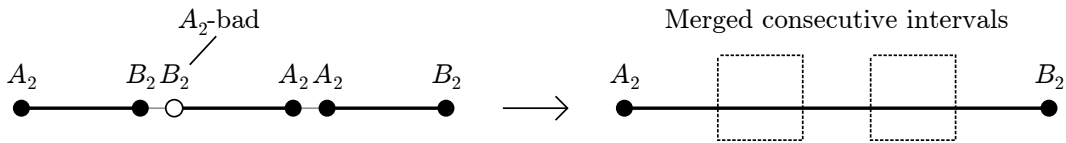}
  \caption{\footnotesize Original edges of $P$ are used to merge
    consecutive intervals in the same block, so that a bad endpoint
    can be eliminated. }
  \label{fig:path-merge}
\end{figure}

\begin{description}
\item[Step 1.] For every vertex $v \in P \setminus B_1$, expose the
  value of $d_1^+(v; B_1)$.  We show that \whp, the initial and final
  endpoints of $P$ are $B_1$-good, and at most $n \cdot e^{-\sqrt{\log
      n}}$ vertices of $P$ are $B_1$-bad.

\item[Step 2.] Absorb all remaining vertices and covering paths into
  the system of segments, using $G_3$-edges that are aligned with
  $B_1$-good vertices.  (See Figure \ref{fig:path-absorb}.)  This
  removes some segment endpoints, while adding other new endpoints.
  Let $A_2 \cup B_2$ be the new partition of endpoints.  Crucially,
  $B_2 = B_1$, while $|A_2| = |A_1|$ by losing up to
  $\frac{2n}{\sqrt[3]{\log n}}$ vertices, and then adding back the
  same number.  Importantly, every new vertex in $A_2 \setminus A_1$
  will be $B_1$-good.

\item[Step 3.] The system of segments can be grouped into several
  blocks of consecutive segments, in the sense that between
  successive segments in the same block, there is an original edge of
  $P$.  (See Figure \ref{fig:path-strings}.)  Also, the initial and
  final endpoints of each block are always of type $A$, and are all
  $B_1$-good.

\item[Step 4.] For every vertex $b \in B_2$, expose the value of
  $d_1^+(b; A_2)$.  We show that \whp, at most $n \cdot e^{-\sqrt{\log
      n}}$ vertices of $B_2$ are $A_2$-bad.

\item[Step 5.] For each consecutive pair of segments along the same
  block (from Step 3) which has either an $A_2$-endpoint which is
  $B_1$-bad or a $B_2$-endpoint which is $A_2$-bad, merge them,
  together with a neighboring segment in order to maintain parity
  between $A$'s and $B$'s.  (See Figure \ref{fig:path-merge}.)  Let
  $A_3 \cup B_3$ be the final partition of segment endpoints after the
  merging.  We show that \whp, all $a \in A_3$ have $d_1^+(a; B_3)
  \geq \frac{\e_1 \th_1}{200L} \log n$, and all $b \in B_3$ have
  $d_1^+(b; A_3) \geq \frac{\e_1 \th_1}{200L} \log n$.  Furthermore,
  if we were to expose the $D_1$-edges between $A_3$ and $B_3$, then
  each vertex $a \in A_3$ would independently sample $d_1^+(a; B_3)$
  uniformly random neighbors in $B_3$, and similarly for $b \in B_3$.
\end{description}

\noindent This will complete the proof because the final number of
segments is $|A_3| \leq |A_2| = |A_1| \leq \frac{n}{L}$.

\subsubsection{Step 1}

By construction, $|B_1| = (1 - o(1)) \frac{n}{L}$.  Now consider an
arbitrary vertex $v \in P \setminus B_1$.  We have only exposed the
numeric value of $d_1^*(v)$ thus far in our construction, and not
where the $D_1$-out-neighbors are.  So let us now expose the numeric
value of $d_1^+(v; B_1)$, but again, not precisely where the
out-endpoints are.  As we observed in the beginning of Section
\ref{sec:long-path}, we have $d_1^*(v) \geq \frac{\e_1 \th_1}{20} \log
n - 3$.  Our work in the previous section consumes up to one
$D_1$-out-edge at each vertex $v \in P$, when we reveal whether
$\rt{vw} \in D_1^\circ$ in the third step of the algorithm.
Therefore, $d_1^+(v; B_1)$ stochastically dominates $\bin{ \frac{\e_1
    \th_1}{20} \log n - 4, (1 - o(1)) \frac{1}{L} }$.  Hence we can
use Lemma \ref{lem:bin9} to bound that the probability that $d_1^+(v;
B_1)$ is too small.
\begin{align*}
  \pr{ d_1^+(v; B_1) < \frac{1}{9} \cdot \frac{1}{L} \cdot \frac{\e_1 \th_1}{20} \log n }
  &<
  e^{-0.533 \cdot \frac{1}{L} \cdot \frac{\e_1 \th_1}{20} \log n } \\
  \pr{ d_1^+(v; B_1) < \frac{\e_1 \th_1}{180 L} \log n }
  &=
  o( e^{-\frac{\e_1 \th_1}{40 L} \log n } ) \\
  &=
  o( e^{-\sqrt{\log n}} ),
\end{align*}
since $\e_1, \th_1 = \Omega \big( \frac{1}{\sqrt{\log \log n}} \big)$
and $L < \sqrt[6]{\log n}$.  The expected number of such vertices in
$P$ is at most $n$ times this probability.  Applying Markov's
inequality, we conclude that \whp, the number of $B_1$-bad vertices in
$P$ is at most $n \cdot e^{-\sqrt{\log n}}$.  This also shows that the
initial and final endpoints of $P$ are $B_1$-good \whp.

\subsubsection{Step 2}

At this point, our entire vertex set is partitioned as follows.  We
have a collection of rainbow intervals $I_1, \ldots, I_r$, each of
length exactly $L$.  These already consume at least $n -
\frac{n}{\sqrt[3]{\log n}}$ vertices.  Since we discarded the
remainder of $P$, as well as possibly the final interval, we have $r
\geq \frac{1}{L} \big( n - \frac{n}{\sqrt[3]{\log n}} \big) - 2$.  A
separate collection of rainbow paths $Q_1, \ldots, Q_s$ covers all
vertices of $S$.  There are also some remaining vertices.  In this
section, we will use $G_3$-edges to absorb the latter two classes into
the rainbow intervals.

Since we will not use any further $G_2$-edges, but edges from
$D_2^\circ$ take precedence over those from $D_3$, we also now expose
all edges in $D_2^\circ$.  Lemma \ref{lem:max-degree} ensures that
\whp, no vertex is incident to more than $5 \log n$ edges of
$D_2^\circ$.  Condition on this outcome.  Note that by construction,
we have not exposed the locations of any $D_3^\circ$-edges between
vertices outside $S$, although vertices outside $S$ may have up to
three exposed $D_3$-neighbors located in $S$.

We now use a simple greedy algorithm to absorb all residual paths and
vertices into our collection of intervals.  In each step, we find a
pair of $G_3$-edges linking either a new $Q_i$ or a new vertex to two
distinct intervals $I_x$ and $I_y$, using two new colors from $C_3$.
We will ensure that throughout the process, all intervals $I_w$ used
in this way are separated by at least one full interval $I_z$ along
$P$.

The specific procedure is as follows.  Suppose we have already linked
in $t$ paths or vertices, and are considering the next path or vertex
to link in.  Suppose it is a path $Q_i$ (the vertex case can be
treated in an analogous way).  Let $u, v$ be the endpoints of $Q_i$.
We need to find vertices $x,y$ in distinct intervals $I_x$ and $I_y$
such that (i) according to $P$, $I_x$ and $I_y$ are separated by at
least one full interval from each other, and from all $I_z$ previously
used in this stage, (ii) $x$ and $y$ are separated from the endpoints
of $I_x$ and $I_y$ by at least two edges of $P$, and (iii) if $x' \in
I_x$ is the vertex adjacent to $x$ in the direction of the
$B_1$-endpoint of $I_x$, and $y' \in I_y$ is the vertex adjacent to
$y$ in the direction of the $B_1$-endpoint of $I_y$, then both $x'$
and $y'$ are $B_1$-good.  We choose the direction of the
$B_1$-endpoint because $x'$ and $y'$ will become the new
$A_2$-endpoints of shortened intervals; see Figure
\ref{fig:path-absorb} for an illustration.

So, let $F$ be the set of vertices in the intervals which fail
properties (ii) or (iii).  By Step 1, the dominant term arises from
the endpoints because $n \cdot e^{-\sqrt{\log n}} \ll
\frac{n}{\sqrt[6]{\log n}} < \frac{n}{L}$, so $|F| < \frac{5n}{L}$.
Also let $T$ be the set of vertices contained in intervals that are at
least one full interval away from any intervals which have previously
been touched by this algorithm.  Since we observed at the beginning of
this section that the total number of intervals was at least
$\frac{1}{L} \big( n - \frac{n}{\sqrt[3]{\log n}} \big) - 2$, we have
\begin{displaymath}
  |T|
  \geq
  \left(
    \frac{1}{L} \left( n - \frac{n}{\sqrt[3]{\log n}} \right) - 2 - 3\cdot (2t)
  \right) L
  \geq
  n - \frac{8Ln}{\sqrt[3]{\log n}}
  \geq
  n - \frac{8n}{L},
\end{displaymath}
since $t \leq \frac{n}{\sqrt[3]{\log n}}$ and $L < \sqrt[6]{\log n}$.
Let us now find a newly-colored $G_3$-edge from $u$ to a vertex of $T
\setminus F$.  Note that the number of vertices outside $T \setminus
F$ is at most $\frac{13n}{L}$.

We have not yet exposed the specific locations of the $D_3$-neighbors
of $u$, but only know (since $u \not \in S$) that $d_3^+(u) \geq
\frac{\e_3 \th_3}{20} \log n$, and up to three of those
$D_3$-out-neighbors lie within $S$.  Consider what happens when we
expose the location $w$ of one of $u$'s $D_3$-out-neighbors which is
outside $S$.  This will produce a useful $G_3$-edge $uw$ if (i) $w$
lands in $T \setminus B$, (ii) neither $\rt{uw}$ or $\rt{wu}$ appeared
in $D_1^\circ$ or $D_2^\circ$, (iii) $\rt{wu}$ does not appear in
$D_3^\circ$, and (iv) the color of the edge is new.

Let us bound the probability that $w$ fails any of these properties.
We may consider (i)-(iii) together, since we showed that at most
$\frac{13n}{L}$ vertices were outside $T \setminus B$, and we
conditioned on $u$ being incident to at most $5 \log n$ edges of
$D_2^\circ$.  After sampling the location of $w$, we expose whether
$\rt{uw}$ or $\rt{wu}$ appear in $D_1^\circ$; by the same argument as
used in Lemma \ref{lem:algo-red-whp}, the probability of each is at
most $\frac{7 \log n}{n} = o\big( \frac{1}{L} \big)$.  So, the
probability of failing either (i) or (ii) is at most $\frac{14}{L}$.
We conditioned at the beginning of Section \ref{sec:long-path} on
$d_3^+(w) \leq 5 \log n$, so when we expose whether $\rt{wu}$ is in
$D_3^\circ$, we again fail only with probability at most $\frac{7 \log
  n}{n} = o(\frac{1}{L})$.  Finally, when we expose the color of the
new edge, we know that it will be in $C_3$, so the probability that it
is a previously used color is at most
\begin{displaymath}
  (2t) \cdot \frac{1}{\th_3 n}
  <
  \frac{2n}{\sqrt[3]{\log n}}  \cdot \frac{1}{\th_3 n}
  =
  \frac{2}{\th_3 \sqrt[3]{\log n}}
  =
  o \left( \frac{1}{L} \right).
\end{displaymath}
Therefore, the probability that all of the $\geq \frac{\e_3 \th_3}{20}
\log n - 3$ $D_3$-out-edges of $u$ fail is at most
\begin{displaymath}
  \bfrac{15}{L}^{\frac{\e_3 \th_3}{20} \log n - 3}
  =
  \bfrac{L}{15}^3 n^{ - \frac{\e_3 \th_3}{20} \log \frac{L}{15} }
  \leq
  \bfrac{L}{15}^3 n^{-2},
\end{displaymath}
where we used the definition of $L$ in \eqref{eq:L} for the final
bound.  A similar calculation works for $v$, and for the separate case
when we incorporate a new vertex into the intervals.  Therefore,
taking a union bound over the $o(n)$ iterations in linking vertices
and paths, we conclude that our procedure completes successfully \whp.
\proofend

\subsubsection{Step 3}

Step 1 established that \whp, the initial and final vertices of the
long path $P$ are both $B_1$-good, so the original system of segments
can be arranged as a single block, with successive intervals linked by
edges of $P$.  We now prove by induction that after the absorption of
each path or vertex in Step 2, the collection of segments can be
re-partitioned into blocks of segments, such that within each block,
consecutive segments have their endpoints linked via $P$, and the
initial and final endpoints in each block are of type-$A$, and
$B_1$-good.

There are two cases, depending on whether the absorption involves two
segments in the same block (as in Figure \ref{fig:path-absorb}), or in
different blocks (as in Figure \ref{fig:path-strings}).  If the
segments are in the same block, then we can easily divide that block
into two blocks satisfying the condition.  Indeed, in Figure
\ref{fig:path-absorb}, one of the new blocks is the string of segments
between the vertices indicated by $A_2$ in the diagram, and the other
new block is the complement.  This works because within each of the
two new blocks, every edge between successive segments was an edge
between successive segments of the original block, hence in $P$.
Also, of the four initial/final endpoints among the two new blocks,
two of them were the initial/final endpoints of the original block,
and the other two were identified as $B_1$-good vertices, now in
$A_2$.  Therefore, the new block partition satisfies the requirements.

On the other hand, if the absorption involves two segments from
different blocks, then one can re-partition the two blocks into three
new blocks, as illustrated in Figure \ref{fig:path-strings}.  A
similar analysis to above then completes the argument.

\subsubsection{Step 4}

We have not yet revealed anything about the $D_1$-out-neighbors of any
vertices in $B_2 = B_1$; the only thing we know is that they had
$d_1^+ \geq \frac{\e_1 \th_1}{20} \log n$.  For each vertex $b \in
B_1$, let us now expose the numeric value of $d_1^+(b; A_2)$, but
again, not precisely where the out-endpoints are.  Since our
absorption procedure maintained $|A_2| = |A_1| = (1 - o(1))
\frac{n}{L}$, the same argument that we used for Step 1 now
establishes Step 4.

\subsubsection{Step 5}

By Steps 1 and 4, the total number of merges which occur in Step 5 is
at most $O(n \cdot e^{-\sqrt{\log n}})$.  Since $B_3 \subset B_2 =
B_1$ and $A_3 \subset A_2$, for every $a \in A_3$ we can independently
sample $d_1^+(a; B_3)$ using only the value of $d_1^+(a; B_1)$.
Indeed, since $a \in A_3$ was $B_1$-good, it had $d_1^+(a; B_1) \geq
\frac{\e_1 \th_1}{180 L} \log n$.  We will only have $d_1^+(a; B_3) <
\frac{\e_1 \th_1}{200 L} \log n$ if at least $\frac{\e_1 \th_1}{180 L}
\log n - \frac{\e_1 \th_1}{200 L} \log n = \frac{\e_1 \th_1}{1800 L}
\log n$ of those out-neighbors land in $B_1 \setminus B_3$ as opposed
to $B_3$.  Let
\begin{displaymath}
  q = \frac{|B_1 \setminus B_3|}{|B_1|}
  =
  O \left(
    \frac{n \cdot e^{-\sqrt{\log n}}}{n/L}
  \right)
  =
  e^{-(1 - o(1))\sqrt{\log n}}.
\end{displaymath}
Note that
\begin{align*}
  \pr{
    \bin{\frac{\e_1 \th_1}{180 L} \log n, q}
    \geq
    \frac{\e_1 \th_1}{1800 L} \log n
  }
  &\leq
  \binom{ \frac{\e_1 \th_1}{180 L} \log n}
{\frac{\e_1 \th_1}{1800 L} \log n}
  q^{\frac{\e_1 \th_1}{1800 L} \log n} \\
  &\leq
  \left(
    10eq
  \right)^{\frac{\e_1 \th_1}{1800 L} \log n} \\
  &=
  e^{- (1-o(1)) \sqrt{\log n} \cdot \frac{\e_1 \th_1}{1800 L} \log n} \\
  &<
  e^{- \Omega \left( (\log n)^{4/3} / \log \log n \right)} \\
  &=
  o(n^{-1}).
\end{align*}
Therefore, a final union bound establishes that \whp, every $a \in
A_3$ has $d_1^+(a; B_3) \geq \frac{\e_1 \th_1}{200 L} \log n$.

A similar argument establishes the bound for $d_1^+(b; A_3)$, for $b
\in B_3$, because we had $A_3 \subset A_2$, and had only previously
exposed the value of $d_1^+(b; A_2)$.  The last claim in Step 5 is
clear from our order of exposure.

\subsection{Proof of Lemma \ref{lem7}}

In this final stage of the proof, we use $G_1$ to link together the
endpoints of the system of segments $I_1, \ldots, I_r$ created by
Lemma \ref{lem5}.  As described in the overview (Section
\ref{sec:main-steps}), we construct an auxiliary directed graph
$\Gamma$.  Importantly, no $D_1^\circ$-edges have been revealed
between the endpoint sets $A$ and $B$, so we may now specify a model
for the random $r$-vertex digraph $\Gamma$.

Indeed, consider a vertex $w_k \in \Gamma$, $1 \leq k \leq r$, and let
$a, b$ be the $A$- and $B$-endpoints of the corresponding interval
$I_k$.  We first generate a set $E_1$ of $\Gamma$-edges by sending
exactly $d_1^+(b; A)$ directed edges out of $w_k$, and exactly
$d_1^+(a; B)$ directed edges in to $w_k$.  This is analogous to the
$d$-in, $d$-out model, except that not all degrees are equal.  Some
directed edges will be generated twice; let $F_1$ be that subset, but
keep only one copy in $E_1$.  Color every edge of $E_1$ independently
from $C_1$.  Finally, generate a random subset $F_2 \subset E_1
\setminus F_1$ by independently sampling each edge of $E_1 \setminus
F_1$ with probability $\frac{1}{2} p_1 \cdot \big( \frac{\th_1}{ 1 +
  \th_1 + \th_2 + \th_3} \big)$.  Let $E_1 \setminus F_2$ be the final
edge set of $\Gamma$.

The reason for the removal of $F_2$ is that some of the
initially-generated edges of $\Gamma$ will find conflicts once
$D_1^\circ$ is generated.  Indeed, every edge $\rt{w_j w_k}$ that we
have placed in $\Gamma$ corresponds to an edge $\rt{ba} \in D_1$ or
$\rt{ab} \in D_1$ (or both), for some $b \in B$, $a \in A$.  When both
do not occur, and only $\rt{ba} \in D_1$ (say), then we need to expose
whether $\rt{ab} \in D_1^\circ \setminus D_1$; if it is in $D_1^\circ
\setminus D_1$, it removes $\rt{w_j w_k}$ from $\Gamma$ with
probability $1/2$.

To simplify notation, let $\delta^+(w_k)$ and $\delta^-(w_k)$ be the
numbers of out- and in-edges that are generated at $w_k$ to build the
initial edge set $E_1$.  They correspond to $d_1^+(b; A)$ and
$d_1^+(a; B)$ above, and have therefore been revealed by our previous
exposures.  Importantly, we have the bounds
\begin{displaymath}
  \frac{\e_1 \th_1}{200 L} \log n
  \leq
  \delta^\pm(w_k)
  \leq
  5 \log n.
\end{displaymath}

It is more convenient to restrict our attention to a smaller subset
$E_2 \subset E_1$ which is itself already rainbow; then, every
ordinary directed Hamilton cycle will automatically be rainbow.  For
this, we expose at every vertex $w_k$ what the colors of the
$\delta^+(w_k)$ out-edges and $\delta^-(w_k)$ in-edges will be, but
not their locations.

\begin{lemma}
  \label{lem:select-3-colors}
  Suppose that $\delta^\pm(w_k) \geq \frac{\e_1 \th_1}{200L} \log n$
  for all $1 \leq k \leq r$.  Then \whp, it is possible to select 3
  out-edges and 3 in-edges from each $w_k$ so that all $6r$ selected
  colors are distinct.
\end{lemma}

\proofstart Construct an auxiliary bipartite graph $H$ with vertex
partition $W \cup C_1$, where $W = \{w_1^+, w_1^-, \ldots, w_r^+,
w_r^-\}$.  Place an edge between $w_k^+$ and $c$ if one of $w_k$'s
$\delta^+(w_k)$ out-edges has color $c$.  Edges between $w_k^-$ and
$c$ are defined with respect to $w_k$'s in-edge colors.  The desired
result is a perfect 1-to-3 matching in $H$.  For this, we apply the
1-to-3 version of Hall's theorem: we must show that for every $X
\subset W$, we have $|N(X)| \geq 3|X|$, where $N(X)$ is the union of
the $H$-neighborhoods of all vertices of $X$.

This follows from a standard union bound.  Indeed, fix an integer $1
\leq x \leq 2r$, and consider an arbitrary pair of subsets $X \subset
W$ and $Y \subset C$, with $|X| = x$ and $|Y| = 3x$.  The probability
that $N(X) \subset Y$ (in $H$) is at most
\begin{displaymath}
  \left[
    \bfrac{3x}{\th_1 n}^{\frac{\e_1 \th_1}{200 L} \log n}
  \right]^x
  \, .
\end{displaymath}
The innermost term is the probability that a random color from $C_1$
is in $Y$.  The exponents come from the fact that each vertex of $w_k
\in X$ samples at least $\delta^\pm(w_k) \geq \frac{\e_1 \th_1}{200L}
\log n$ colors for its out- and in-neighbors.

Multiplying this bound by the number of ways there are to select $X$
and $Y$, and using $r \leq \frac{n}{L}$, we find that the probability
of failure for a fixed $x$ is at most
\begin{align*}
  \binom{2r}{x}
  \binom{\th_1 n}{3x}
  \left[
    \bfrac{3x}{\th_1 n}^{\frac{\e_1 \th_1}{200 L} \log n}
  \right]^x
  &\leq
  \bfrac{2er}{x}^x
  \bfrac{\th_1 n}{3x}^{3x}
  \left[
    \bfrac{3x}{\th_1 n}^{\frac{\e_1 \th_1}{200 L} \log n}
  \right]^x \\
  &=
  \left[
    \bfrac{2en}{Lx}
    \bfrac{\th_1 n}{3x}^3
    \bfrac{3x}{\th_1 n}^{\frac{\e_1 \th_1}{200 L} \log n}
  \right]^x
  \, .
\end{align*}
We will sum this over all $1 \leq x \leq 2r$.  The outer exponent
allows us to bound this by a decreasing geometric series, so it
suffices to show that the interior of the square bracket is uniformly
$o(1)$ for all $1 \leq x \leq 2r$.  Indeed, observe that the exponent
of $x$ inside the bracket is $\frac{\e_1 \th_1}{200 L} \log n - 4 >
0$, so it is maximized at $x = 2r \leq \frac{2n}{L}$.  Yet
\begin{align*}
  \bfrac{2en}{L(2n/L)}
  \bfrac{\th_1 n}{3(2n/L)}^3
  \bfrac{3(2n/L)}{\th_1 n}^{\frac{\e_1 \th_1}{200 L} \log n}
  &=
  (e)
  \bfrac{6}{\th_1 L}^{\frac{\e_1 \th_1}{200 L} \log n - 3}
\end{align*}
Since \eqref{eq:L} ensures that $L \geq \frac{7}{\th_1}$, we have
$\frac{6}{\th_1 L} \leq \frac{6}{7}$.  Yet the exponent $\frac{\e_1
  \th_1}{200 L} \log n - 3$ tends to infinity as $n$ grows, so we
indeed obtain a uniform upper bound of $o(1)$ for all $1 \leq x \leq
2r$.  Therefore, \whp, every subset $X \subset W$ has $|N(X)| > 3|X|$,
and the 1-to-3 version of Hall's theorem establishes the desired
result.  \proofend

\vspace{3mm}

Recall from the beginning of this section that the final edge set of
$\Gamma$ is $E_1 \setminus F_2$.  Let $E_2$ be the set of $6r$ edges
selected by Lemma \ref{lem:select-3-colors}.  This corresponds to a
copy of $D_{\text{3-in},\text{3-out}}$.  Unfortunately, in our model
we still need to expose the locations of the remaining
$\delta^\pm(w_k) - 3$ remaining in- and out-edges at every vertex
$w_k$.  It is possible that an edge of $E_2$ may be generated again in
this stage.  That edge would then have $1/2$ probability of receiving
the color of the new copy, which would not be in our specially
constructed rainbow set.  To account for this, let $E_3$ be the set of
edges generated by exposing these remaining in- and out-edges, so that
$E_2 \cup E_3 = E_1$.  It suffices to find a directed Hamilton cycle
in $E_2 \setminus (E_3 \cup F_2)$.

It is not convenient to work directly with $E_3$ or $F_2$, because
they depend on the result of $E_2$.  Let $F_3$ be the random directed
graph $D_{r,q}$ defined by sampling each edge with probability $q =
\frac{130L \log n}{n}$, and generated independently of $E_2$.
Fortunately, we can control $F_3$ instead.

\begin{lemma}
  \label{lem:coupling}
  There is a coupling of the probability space such that $E_3 \cup F_2
  \subset F_3$ \whp.
\end{lemma}

\proofstart The set $F_2$ is a random subset of $E_1 \setminus F_1$
obtained by independently sampling each edge with a probability of at
most $\frac{\log n}{n}$, so it is clearly contained in a copy of
$D_{r,\frac{\log n}{n}}$ that is generated independently of $E_2$.
Next, although $E_3$ is exposed after $E_2$, it is still contained in
a copy of $D_{(6 \log n)\text{-in},(6 \log n)\text{-out}}$ that is
generated independently of $E_2$.  Indeed, at a vertex $w_k$, we
generate $E_3$ by exposing $\delta^\pm(w_k) - 3$ new out- and
in-edges, but $\delta^\pm(w_k) \leq 5 \log n$.

We observe a standard coupling which realizes $D_{(6 \log
  n)\text{-in},(6 \log n)\text{-out}}$ as a subgraph of
$D_{r,\frac{120 \log n}{r}}$ \whp.  For this, consider the following
system for generating a random directed graph.  For every ordered pair
of vertices $(u, v)$, generate two independent Bernoulli random
variables $I_{u,v}^+$ and $I_{u,v}^-$, each with probability parameter
$\frac{60 \log n}{r}$.  Create the directed edge $\rt{uv}$, if and
only if at least one of $I_{u,v}^+$ or $I_{v,u}^-$ took the value 1.
This is clearly contained in $D_{r,q}$.  For each vertex $u$, let
$D^+(u)$ be the number of other vertices $v$ for which $I_{u,v}^+ =
1$, and let $D^-(u)$ be the number of $I_{u,v}^- = 1$.  These are all
distributed as $\bin{r-1, \frac{60 \log n}{r}}$.  Lemma \ref{lem:bin9}
establishes that for fixed $u$, the probability that $D^+(u) < 6 \log
n$ is at most $e^{- 0.533 \cdot (60 - o(1)) \log n}$, and similarly
for $D^-(u)$.  Therefore, a union bound establishes that \whp, all
$D^+(u), D^-(u) \geq 6 \log n$.  Conditioning on the values of
$D^+(u)$ and $D^-(u)$, we see that when the indicators are revealed,
this indeed contains $D_{(6 \log n)\text{-in},(6 \log n)\text{-out}}$.
The result then follows by recalling that $r = (1-o(1)) \frac{n}{L}$,
and taking the union of $D_{r,\frac{120L \log n}{n}}$ with another
independent $D_{r,\frac{\log n}{n}}$ to cover $F_2$.  \proofend

\vspace{3mm}

We are now ready to finish the final lemma in our proof of Theorem
\ref{th1}.

\vspace{3mm}

\noindent \textbf{Proof of Lemma \ref{lem3}.}\, We have established
that it suffices to find an ordinary directed Hamilton cycle in $E_2
\setminus F_3$, without regard to color.  Conveniently, $F_3$ is now
independent of $E_2 \sim D_{\text{3-in},\text{3-out}}$.  This
independence allows us to first expose how many of the 3-in, 3-out
edges at each vertex will be in $F_3$, and then only expose the
locations of those that are not in $F_3$.  At each vertex $u$, the
probability that more than one of the 6 in- or out-edges is in $F_3$
is at most
\begin{displaymath}
  \binom{6}{2}
  \left(
    \frac{130L \log n}{n}
  \right)^2
  =
  o(n^{-1}),
\end{displaymath}
so a union bound establishes that \whp, $E_2 \setminus F_3$ contains a
copy of $D_{\text{2-in},\text{2-out}}$.  This is known to be
Hamiltonian \whp\ by Theorem \ref{thm:2-in-2-out} from Section
\ref{sec:main-steps}, so our proof is complete.  \proofend

\section{Concluding remarks}

Our main contribution, part (b) of Theorem \ref{th1}, sharpens the
earlier result of Cooper and Frieze \cite{CF} to achieve optimal
first-order asymptotics.  As we mentioned in the introduction, we
suspect that our result can be further sharpened within the $o(1)$
term.  We do not push to optimize our bounds on $\e, \th$ because it
is not clear that incremental improvements upon our current approach
will be sufficiently interesting.  Instead, we would be more
interested in determining whether one can extend the celebrated
``hitting time'' result of Bollob\'as \cite{Bol-hitting-time} , which
states that one can typically find a Hamilton cycle in the random
graph process as soon as the minimum degree reaches two.

\begin{question}
  Consider the edge-colored random graph process in which
  $e_1,e_2,\ldots,e_N,\,N=\binom{n}{2}$ is a random permutation of the
  edges of $K_n$. The graph $G_m$ is defined as
  $([n],\set{e_1,e_2,\ldots,e_m})$.  Let each edge receive a random
  color from a set $C$ of size at least $n$. Then \whp, the first time
  that a rainbow Hamilton cycle appears is precisely the same as the
  first time that the minimum degree of $G_m$ is at least two
  \textbf{and} at least $n$ colors have appeared.
\end{question}

Although this may be out of reach at the moment, another natural
challenge is to settle the problem in either the case when the number
of edges is just sufficient for an ordinary Hamilton cycle, or in the
case when the number of colors is, say, exactly $n$.  Part (a) of our
theorem answered the latter question when $n$ was even.  We believe
that it is probably also true when $n$ is odd.  It would be nice to
prove that extension, either directly or by reducing the divisibility
condition in Theorem \ref{th:loose-3}.


\begin{thebibliography}{99}

\bibitem{AFR} M. Albert, A.M. Frieze and B. Reed, Multicoloured
  Hamilton cycles, \emph{Electronic Journal of Combinatorics}
  \textbf{2} (1995), R10.

\bibitem{AS} N. Alon and J. Spencer, \textbf{The Probabilistic
    Method}, 3rd ed., Wiley, New York, 2007.

\bibitem{Bolbook} B. Bollob\'as, \textbf{Random Graphs}, 2nd ed.,
  Cambridge University Press, 2001.

\bibitem{Bol-hitting-time} B. Bollob\'as, The evolution of sparse
  graphs, in: \textbf{Graph Theory and Combinatorics}, Proc. Cambridge
  Combinatorial Conf. in honour of Paul Erd\H{o}s (B. Bollob\'as,
  ed.), Academic Press, 1984, 35--57.

\bibitem{BolFrieze-hitting-time} B. Bollob\'as and A. Frieze, On matchings and
  Hamiltonian cycles in random graphs, in: \textbf{Random Graphs '83
    (Poznan, 1983)}, North-Holland Math. Stud., 118, North-Holland,
  Amsterdam (1985), 23--46.

\bibitem{CF-2-inout} C. Cooper and A.M. Frieze, Hamilton cycles in
  random graphs and directed graphs, \emph{Random Structures and
    Algorithms} \textbf{16} (2000), 369--401.

\bibitem{CF1} C. Cooper and A.M. Frieze, Multicoloured Hamilton cycles
  in random graphs: an anti-Ramsey threshold, \emph{Electronic Journal
    of Combinatorics} \textbf{2} (1995), R19.

\bibitem{CF} C. Cooper and A.M. Frieze, Multi-coloured Hamilton cycles
  in random edge-coloured graphs, {\em Combinatorics, Probability and
    Computing} \textbf{11} (2002), 129--134.

\bibitem{DudekFr} A. Dudek and A.M. Frieze, Loose Hamilton cycles in
  random $k$-uniform hypergraphs, submitted.

\bibitem{FrH} A.M. Frieze, Loose Hamilton cycles in random 3-uniform
  hypergraphs, \emph{Electronic Journal of Combinatorics} \textbf{17}
  (2010), N28.

\bibitem{FK-pack} A.M. Frieze and M. Krivelevich, Packing Hamilton
  cycles in random and pseudo-random hypergraphs, submitted.

\bibitem{FKL-pack} A.M. Frieze, M. Krivelevich, and P. Loh, Packing
  tight Hamilton cycles in 3-uniform hypergraphs, \emph{Random
    Structures and Algorithms}, to appear.

\bibitem{ENR} P. Erd\H{o}s, J. Ne\v{s}et\v{r}il, and V. R\"odl, Some
  problems related to partitions of edges of a graph, in:
  \textbf{Graphs and Other Combinatorial Topics}, Teubner, Leipzig
  (1983), 54--63.
  
\bibitem{DFrieze} A.M. Frieze, An algorithm for finding Hamilton
  cycles in random digraphs, \emph{Journal of Algorithms} \textbf{9}
  (1988), 181--204.

\bibitem{FM} A.M. Frieze and B. D. McKay, Multicoloured trees in
  random graphs, {\em Random Structures and Algorithms} \textbf{5}
  (1994), 45--56.


\bibitem{JW} S. Janson and N. Wormald, Rainbow Hamilton cycles in
  random regular graphs, {\em Random Structures Algorithms}
  \textbf{30} (2007), 35--49.

\bibitem{HT} G. Hahn and C. Thomassen, Path and cycle sub-Ramsey
  numbers, and an edge colouring conjecture, \emph{Discrete
    Mathematics}. \textbf{62} (1986), 29--33.


\bibitem{KS-Hamthresh} J. Koml\'os and E. Szemer\'edi, Limit
  distribution for the existence of Hamiltonian cycles in a random
  graph, \emph{Discrete Mathematics} \textbf{43} (1983), 55--63.

\bibitem{McD-digraphs1} C.J.H. McDiarmid, General percolation and
  random graphs, \emph{Adv. Appl. Probab.} \textbf{13} (1981), 40--60.

\bibitem{McD-digraphs2} C.J.H. McDiarmid, General first-passage
  percolation, \emph{Adv. Appl. Probab.} \textbf{15} (1983), 149--161.

\bibitem{MW} B.D. McKay and N.C. Wormald, Asymptotic enumeration by
  degree sequence of graphs with degree $o(n^{1/2})$,
  \emph{Combinatorica} \textbf{11} (1991), 369--382.

\bibitem{DLV} W. Fernandez de la Vega, Long paths in random graphs, {\em
    Studia Math. Sci. Hungar.} (1979), 335--340.

\bibitem{RW-3reg} R.W. Robinson and N.C. Wormald, Almost all regular
  graphs are Hamiltonian, \emph{Random Structures and Algorithms}
  \textbf{5} (1994), 363--374.

\bibitem{Rue} R. Rue, Comment on \cite{AFR}.
\emph{Electronic Journal of Combinatorics}
  \textbf{2} (1995).

\end{thebibliography}
\end{document}